\documentclass[11pt,a4paper]{amsart}
\usepackage{amssymb}
\usepackage{amscd}
\usepackage{amsmath}\usepackage{amsfonts}
\usepackage{pstricks,pst-node}
\usepackage{mathrsfs}
\usepackage{cite}

\usepackage{mathrsfs}

\numberwithin{equation}{section}
\swapnumbers
\theoremstyle{plain}
\newtheorem{theorem}[equation]{Theorem}
\newtheorem{corollary}[equation]{Corollary}
\newtheorem{lemma}[equation]{Lemma}

\newtheorem{proposition}[equation]{Proposition}

\newtheorem*{thm}{The Carter--Payne Theorem}

\newtheorem{case}{Case}

\theoremstyle{definition}
\newtheorem{definition}[equation]{Definition}

\theoremstyle{remark}

\makeatletter
\def\enumerate{\begingroup\ifnum\@enumdepth>3\@toodeep\else
      \advance\@enumdepth\@ne
      \edef\@enumctr{enum\romannumeral\the\@enumdepth}%
      \topsep\z@\parskip\z@
      \list{\csname label\@enumctr\endcsname}
        {\@nmbrlisttrue\let\@listctr\@enumctr
         \parsep\z@\itemsep\z@\topsep\z@
         \setcounter{\@enumctr}{0}
         \def\makelabel##1{\hss\llap{\rm ##1}}
       }\fi
}
\makeatother

\def\crulefill{\leavevmode\leaders\hrule height 1pt\hfill\kern 0pt}
\long\def\QUERY#1{%
 \leavevmode\newline\noindent$\star\star\star$%
     \thinspace\textsf{Comment/Query}\crulefill\newline%
     \space #1\leavevmode\newline\hbox to 120mm{\crulefill}%
     $\star\star\star$\newline
}

\let\bar=\overline

\def\N{\mathbb N}
\def\Z{\mathbb Z}

\newcommand{\mf}[1]{\mathfrak{#1}}
\newcommand{\down}[2]{T_{\downarrow}(#1,#2)}
\newcommand{\up}[2]{T_{\uparrow}(#1,#2)}
\newcommand{\h}{\ensuremath{\mathscr{H}}}
\newcommand{\sym}[1]{\ensuremath{\mathfrak{S}_{#1}}}

\newcommand{\Hom}{\operatorname{Hom}}
\newcommand{\Ker}{\operatorname{Ker}}
\newcommand{\Image}{\operatorname{Im}}

\newcommand{\mat}[2]{\genfrac{[}{]}{0pt}{}{#1}{#2}}

\begin{document}\bibliographystyle{andrew}
\title{Some $q$-analogues of the Carter--Payne theorem}

\author{Sin\'{e}ad Lyle}
\address{School of Mathematics and Statistics F07, University of
Sydney, NSW 2006, Australia.}
\email{s.lyle@maths.usyd.edu.au}
\thanks{}
\begin{abstract}
We prove a $q$-analogue of the Carter--Payne theorem for the two special cases corresponding to moving an arbitrary number of nodes between adjacent rows, or moving one node between an arbitrary number of rows. As a consequence, we show that these homomorphism spaces are one dimensional when $q \neq -1$.   We apply these results to complete the classification of the reducible Specht modules for the Hecke algebras of the symmetric groups when $q \neq-1$.  Our methods can also be used to determine certain other pairs of Specht modules between which there is a homomorphism. In particular, we describe the homomorphism space $\Hom_\h(S^{(n)},S^\mu)$ for an arbitrary partition $\mu$.
\end{abstract}
\maketitle

\section{Introduction}
Let $F$ be a field, $q$ an invertible element of $F$ and $n$ a positive integer. We consider the representations of the Hecke algebra $\h=\h_{F,q}(\sym{n})$.  For each partition $\lambda$ of $n$, we define a $\h$-module $S^{\lambda}$, called a Specht module; it is well--known that when $\h$ is semisimple, the modules $\{S^\lambda \mid \lambda \text{ is a partition of } n\}$ form a complete set of pairwise non--isomorphic irreducible $\h$-modules.  It is an important open problem to determine the homomorphism spaces $\Hom_\h(S^\lambda,S^\mu)$, for $\lambda$ and $\mu$ partitions of $n$.
The most famous result of this kind for the symmetric groups (that is, the case $q=1$) is the Carter--Payne theorem. 

\begin{thm} [{\cite{Payne}, p.~425}]\label{CarterPayne}
Let $\h\cong F\sym{n}$, where $F$ is a field of characteristic $p> 0$.  Choose $\gamma >0$ and take $\mu$ and $\lambda$ to be partitions of $n$ such that 
\[\lambda_i = \begin{cases}
\mu_i+\gamma & \text{if } i=a, \\
\mu_i-\gamma & \text{if } i=b, \\
\mu_i & \text{otherwise},
\end{cases}\] for some $a<b$.  
Suppose that $\mu_a - \mu_b +b-a +\gamma\equiv 0 \mod p^{\ell_p(\gamma)}$, where $\ell_p(\gamma)$ is the smallest positive integer such that $p^{\ell_p(\gamma)}>\gamma$.
Then there exists a non--zero $\h$-homomorphism $\hat{\Theta}:S^{\lambda} \rightarrow S^{\mu}$.
\end{thm}
Although widely conjectured, no $q$-analogue of the full Carter--Payne theorem is known. In this paper, we prove such an analogue in two important special cases, namely when $b=a+1$ or when $\gamma=1$.  Combinatorially, this corresponds to moving an arbitrary number of nodes between adjacent rows, or moving one node between an arbitrary number of rows.  

It turns out that our proof has interesting implications.
We will turn our attention to the classification of the reducible Specht modules for the algebras $\h_{F,q}(\sym{n})$ when $q \neq -1$.  Recent work (see Section \ref{irred}) has resulted in the completion of this classification when $q=1$. Combining our main result, Theorem \ref{image}, with previously published work, we prove that a Specht module is reducible if and only if it is `$(e,p)$-reducible', as defined in Definition \ref{epdefn}.

We will also use our discussion of semistandard homomorphisms to obtain a description of the Specht modules which contain a submodule isomorphic to the trivial module $S^{(n)}$.  (See Theorem \ref{maphappens}.)  Although seemingly elementary, we believe that it is the first time that such a description has been given.  The result and the methods used turn out to be exact analogues of the work of James \cite{James}, Theorem 24.4 for the symmetric groups.  
  
Our proof of these Carter--Payne $q$-analogues is constructive; we will write down the maps in question. This is simple when $b=a+1$, but when $\gamma =1$, it turns out that there is also an elegant formula in terms of semistandard homomorphisms, given in Theorem \ref{addthm}.
 When $q=1$, an explicit description of all of the Carter--Payne homomorphisms has been given by Fayers and Martin~\cite{FayersMartin:homs}; and Ellers and Murray~\cite{EllersMurray} (independently of this work) have shown that when $\gamma =1$, the dimension of the homomorphism space is at most one. In the presence of the Carter--Payne theorem, they have therefore been able to write down a map which, when $q=1$, agrees with Theorem \ref{addthm}.
However, homomorphisms between Specht modules for arbitrary Hecke algebras are not well understood, even for the Hecke algebras $\h_{\mathbb{C},q}(\sym{n})$ for which the decomposition matrices can be computed~\cite{ariki,LLT}.  
In Proposition \ref{comb}, we give the first step towards a general method for studying homomorphisms between Specht modules, namely a way of combining certain important homomorphisms.  This is a $q$-analogue of \cite{FayersMartin:homs}, Lemma 5 which was heavily used throughout that paper.    

Using Proposition \ref{comb}, and other combinatorial methods, we may manipulate semistandard homomorphisms. This approach has been widely used to study homomorphims between Specht modules for the symmetric groups.  Even though we shall use only the classic theory of Dipper and James~\cite{DJ:reps}, it is the first time that it has been adapted for arbitrary Hecke algebras of symmetric groups.  

\section{The Hecke algebras}\label{hecke}
We begin with some standard definitions and notation, most of which can be found in~\cite{DJ:reps}.
Let $\sym{n}$ denote the symmetric group on $n$ letters and for $1 \leq i <n$, let $s_i$ denote the basic transposition $(i,i+1)$, so that $\{s_i \mid 1 \leq i <n\}$ generates $\sym{n}$. For a permutation $w \in \sym{n}$, the length $\ell(w)$ of $w$ is defined to be the smallest value of $k$ such that $w=s_{i_1}s_{i_2}\ldots s_{i_k}$ for some basic transpositions $s_{i_j}$; note that for $w \in \sym{n}$ and $1 \leq i <n$ we have 
\begin{align}\label{length}
\ell(s_i w)  &= \begin{cases} \ell(w)+1 & \text{ if } iw < (i+1)w, \\
\ell(w)-1 & \text{ if } iw > (i+1)w,
\end{cases}\\
\ell(w s_i)  &= \begin{cases} \ell(w)+1 & \text{ if } iw^{-1} < (i+1)w^{-1}, \\
\ell(w)-1 & \text{ if } iw^{-1} > (i+1)w^{-1}.
\end{cases}
\end{align}

Let $F$ be a field of characteristic $p \geq 0$ and $q$ an an invertible element of $F$.  Define $e>1$ to be minimal such that $1+q+\ldots+q^{e-1}=0$, with $e=\infty$ if no such integer exists.  We define the Hecke algebra $\h=\h_{F,q}(\sym{n})$ to be the associative $F$-algebra with basis
$\{T_w \mid w\in\sym n\}$ and multiplication determined by
$$
T_{w}T_{s_i}=\begin{cases} T_{ws_i},           & \text{if } 
\ell(ws_i)=\ell(w)+1,\\ 
qT_{ws_i}+(q-1)T_w,  &\text{if } \ell(ws_i)=\ell(w)-1, \end{cases}$$
where $w\in\sym n$ and $1 \leq i <n$.
Then $\h$ is generated by the elements $T_{s_1},T_{s_2},\ldots,T_{s_{n-1}}$.  For convenience, we will often write $T_i$ to denote $T_{s_i}$.

Let $\lambda$ be a composition of $n$.  
The diagram of $\lambda$ is the set of nodes
\[[\lambda] = \{(i,j)\mid 1 \leq i \mbox{ and } 1 \leq j\leq \lambda_i\}.\]
A $\lambda$-tableau consists of $[\lambda]$ with the nodes replaced with integers; unless otherwise specified we assume that the nodes are replaced by the elements of $\{1,2,\dots,n\}$ in some order. It is said to be row standard if its entries increase across the rows. The symmetric group acts on the right on the set of $\lambda$-tableaux by permuting the entries.  Define $\mf{t}^\lambda$ to be the row standard $\lambda$-tableau with $1,2,\ldots,n$ entered in order along its rows, and  $\mf{t}_\lambda$ to be the row standard $\lambda$-tableau with $1,2,\ldots,n$ entered in order down its columns.  Let $w_\lambda$ be the permutation that sends $\mf{t}^\lambda$ to $\mf{t}_\lambda$.  So, if $\lambda = (3,2)$ then
\[ \mf{t}^\lambda = \begin{array}{lll}1&2&3\\4&5&\end{array}, \qquad \mf{t}_\lambda=\begin{array}{lll}1&3&5\\2&4&\end{array}\]
and $w_\lambda=(2,3,5,4)$.
Let $\sym{\lambda}$ denote the row--stabilizer of $\mf{t}^\lambda$.  Hence define 
\begin{align*}
x_\lambda & = \sum_{w \in \sym{\lambda}} T_w, \\
y_\lambda & = \sum_{w \in \sym{\lambda}} (-q)^{-\ell(w)}T_w.
\end{align*}
Note that if $v \in \sym{\lambda}$ then
$x_{\lambda} T_v = q^{\ell(v)} x_\lambda= T_v x_\lambda \text{ and }  y_{\lambda} T_v = (-1)^{\ell(v)} y_\lambda=T_v y_\lambda.$
 We define the permutation module $M^\lambda$ to be the right $\h$-module $x_\lambda \h$.  Suppose now that $\lambda$ is a partition.  Let $\lambda'$ denote the partition conjugate to $\lambda$, that is, the partition obtained by swapping the rows and columns of $[\lambda]$. Hence define the Specht module $S^\lambda$ to be the right $\h$-module $x_\lambda  T_{w_\lambda} y_{\lambda'} \h$.  Clearly $S^\lambda$ is a submodule of $M^\lambda$.

Set $\mathscr{D}_\lambda =\{d \in \sym{n} \mid \mf{t}^\lambda d \text{ is a row standard $\lambda$-tableau} \}$.  Then $\mathscr{D}_\lambda$ is a complete set of right coset representatives of $\sym{\lambda}$ in $\sym{n}$, consisting of the unique element of minimal length from each coset.  The elements $\{x_\lambda T_{d} \mid d \in \mathscr{D}_\lambda \}$ form a basis of $M^\lambda$.  
Note that if $v \in \sym{\lambda}$ and $d \in \mathscr{D}_\lambda$ then $\ell(vd)=\ell(v)+\ell(d)$.  Then if $w \in \sym{n}$, we may write $T_w=T_v T_d$ where $v \in \sym{\lambda}$ and $d \in \mathscr{D}_\lambda$ sends $\mf{t}^\lambda$ to the row standard $\lambda$-tableau obtained by reordering the rows of $\mf{t}^\lambda w$. Thus $x_\lambda T_w = q^{\ell(v)}x_\lambda T_d$.

\begin{lemma}[{\cite{DJ:reps}, Lemma 3.2}] \label{lemkill}
Let $d \in \mathscr{D}_\lambda$ and $1 \leq i <n$.  Then
\[x_\lambda T_d T_{s_i} = \begin{cases}
qx_\lambda T_d & \text{ if $i,i+1$ belong to the same row of $\mf{t}^\lambda d$}, \\
x_\lambda T_{d s_i} &\text{ if the row index of $i$ in  $\mf{t}^\lambda d$ is less than that of $i+1$}, \\
q x_\lambda T_{d s_i} + (q-1)x_\lambda T_d &\text{ otherwise}.
\end{cases}\]
\end{lemma}

Let $\lambda$ and $\mu$ be compositions of $n$.  A $\lambda$-tableau of type $\mu$ is a tableau of shape $\lambda$ with $\mu_i$ entries equal to $i$, for each $i$.  For a tableau $A$ (of arbitrary type and shape) we write $A(a,b)$ for the entry in the $(a,b)$-place of $A$.
The tableau $A$ is said to be row standard if its entries increase along the rows, and semistandard if its entries both increase along the rows and strictly increase down the columns.  
Let $\mathcal{T}(\lambda,\mu)$ denote the set of $\lambda$-tableaux of type $\mu$, and $\mathcal{T}_0(\lambda,\mu)$ denote the set of semistandard $\lambda$-tableaux of type $\mu$.  We define an equivalence relation $\sim_r$ on $\mathcal{T}(\lambda,\mu)$ by saying that $A \sim_r B$ if for all $i$, row $i$ of $A$ contains the same numbers as row $i$ of $B$.  

Let $A \in \mathcal{T}(\lambda,\mu)$.  Define $1_A \in \sym{n}$ to be the permutation obtained by setting $\mf{t}^\mu 1_A$ to be the row standard $\mu$-tableau for which $i$ belongs to row $r$ if the place occupied by $i$ in $\mf{t}^\lambda$ is occupied by $r$ in $A$.  Then $A \mapsto 1_A$ gives a bijection between $\mathcal{T}(\lambda,\mu)$ and $\mathscr{D}_\mu$.

For $A \in \mathcal{T}(\lambda,\mu)$, we now define the homomorphism $\Theta_A:M^\lambda \rightarrow M^\mu$.  For all $h \in \h$, 
\[\Theta_A(x_\lambda h) = \bigg( x_\mu \sum_{A' \sim_r A}T_{1_{A'}}\bigg) h.\]
Then $\{\Theta_A \mid A \in \mathcal{T}(\lambda,\mu) \text{ and is row standard}\}$ form a basis of $\Hom_\h(M^\lambda,M^\mu)$.
The way in which these maps were constructed means that there exists $h' \in \h$ such that $ x_\mu \sum_{A'\sim_r A}T_{1_{A'}}= h' x_\lambda$; for our purposes, it is not necessary to describe $h'$.  More details can be found in~\cite{M:ULect}, 4.5.

Suppose that $\lambda$ is a partition and that $\Theta:M^\lambda \rightarrow M^\mu$.  Let $\hat{\Theta}$ denote the restriction of $\Theta$ to $S^\lambda$.  We will repeatedly use the following theorem.
\begin{theorem}[{~\cite{DJ:qWeyl}, Corollary 8.7}] \label{semistand}
Suppose that $\lambda$ is a partition of $n$ and $\mu$ is a composition of $n$. Then $\{\hat{\Theta}_A \mid A \in \mathcal{T}_0(\lambda,\mu)\}$ is a linearly independent subset of $\Hom_\h(S^\lambda,M^\mu)$, and if either $e \neq 2$ or $\lambda$ is 2-regular (that is, no 2 parts of $\lambda$ are the same length) then $\{\hat{\Theta}_A \mid A \in \mathcal{T}_0(\lambda,\mu)\}$ is a basis of $\Hom_\h(S^\lambda,M^\mu)$.  
\end{theorem}

We now use these homomorphisms to give an alternative description of the Specht module.  Let $\mu$ be a partition.  Take $d$ to be a positive integer and choose $t$ such that $0 \leq t < \mu_{d+1}$.  Let $\nu^{d,t}$ be the composition determined by
\[\nu^{d,t}_i = \begin{cases} \mu_i+\mu_{i+1}-t & \text{if } i=d, \\
t & \text{if } i =d+1, \\
\mu_i & \text{otherwise}. \end{cases} \]
Let $A$ be the row standard $\mu$-tableau of type $\nu^{d,t}$ with all entries in row $i$ equal to $i$, except for row $d+1$ which contains $\mu_{d+1}-t$ entries equal to $d$ and $t$ entries equal to $d+1$.  We write $\psi_{d,t}$ for the homomorphism $\Theta_A:M^\mu \rightarrow M^{\nu^{d,t}}$.  

\begin{theorem}[{~\cite{DJ:reps}, Theorem 7.5}]\label{altdefn}
If $\mu$ is a partition of $n$ then 
\[S^\mu = \bigcap_{d \geq 1} \bigcap_{t=0}^{\mu_{d+1}-1} \Ker \psi_{d,t}. \]
\end{theorem}
We immediately get the following corollary:
\begin{corollary}\label{cor:compose}
Let $\lambda$ and $\mu$ be partitions of $n$ and suppose that $\hat{\Theta}:S^\lambda \rightarrow M^\mu$.  Then $\Image(\hat{\Theta}) \subseteq S^\mu$ if and only if $\psi_{d,t} \hat{\Theta} = 0$ for all $d\geq1$ and $0 \leq t < \mu_{d+1}$.  
\end{corollary}

There are some cases where we may immediately say that $\psi_{d,t} \hat{\Theta} = 0$.

\begin{lemma}\label{nicelemma2}
Suppose that $w \in \sym{n}$ is such that $\mf{t}^\nu w$ contains two entries from the same column of $\mf{t}_\lambda$ in the same row.
Then $x_\nu T_{w} y_{\lambda'} =0$.
\end{lemma}
\begin{proof}
First note that if $s_i \in \sym{\lambda'}$ and $w \in \sym{n}$ then $T_w y_{\lambda'} = - T_{w s_i} y_{\lambda'}$ since if $\ell(w s_i) > \ell(w)$ then
\[ T_w y_{\lambda'} = -T_w T_{s_i} y_{\lambda'} = - T_{w s_i} y_{\lambda'},\] and if  $\ell(w s_i) < \ell(w)$ then 
\[ T_w y_{\lambda'} = T_{(w s_i) s_i} y_{\lambda} = T_ {(w s_i)}T_{s_i} y_{\lambda'} = - T_ {w s_i} y_{\lambda'}.\]
Suppose that $x$ and $y$ lie in the same column of $\mf{t}_\lambda$ and the same row of $\mf{t}^\nu w$ and assume $x<y$.  Let $v = (x+1,x+2,\ldots,y) \in \sym{\lambda'}$.  Then $x_\nu T_w y_{\lambda'} = (-1)^{\ell(v)} x_\nu T_{wv} y_{\lambda'}$. 
  Since $s_x \in \sym{\lambda'}$, we may write $y_{\lambda'} = (I-q^{-1}T_{s_x}) y$ for some $y  \in \sym{\lambda'}$, where $I$ denotes the identity element of $\h$.  The tableau $\mf{t}^\nu wv$ contains the entries $x$ and $x+1$ in the same row.
Now by Lemma \ref{lemkill}, $x_\nu T_{wv} (1-q^{-1}T_{s_x}) =0$.
\end{proof}

\begin{lemma} [{See~\cite{James}, Lemma 3.7}] \label{die2} Define a partial order $\trianglerighteq$ on the set of compositions of $n$ by saying that $\lambda \trianglerighteq \nu$ if and only if 
\[\sum_{i=1}^k \lambda_i \geq \sum_{i=1}^k \nu_i\] for all $k$.
Suppose that $\mf{t}_1$ is a $\lambda$-tableau and $\mf{t}_2$ is a $\nu$-tableau such that for every $i$, the numbers from column $i$ of $\mf{t}_1$ belong to different rows of $\mf{t}_2$.  Then $\lambda  \trianglerighteq \nu$.
\end{lemma}

\begin{proof}
By reordering their parts, we may assume that $\lambda$ and $\nu$ are partitions.  We must place the $\lambda_1'$ numbers from the first column of $\mf{t}_1$ in different rows of $\mf{t}_2$.  Hence $\nu_1' \geq \lambda_1'$.  Next insert the numbers from the second column of $\mf{t}_1$ into differents rows of $\mf{t}_2$.  To do this, we require $\nu_1'+\nu_2' \geq \lambda_1'+\lambda_2'$.  Continuing in this way, we have $\nu' \trianglerighteq \lambda'$.  But it is well known that $\nu' \trianglerighteq \lambda'$ if and only if $\lambda \trianglerighteq \nu$.
\end{proof}

\begin{lemma}\label{die3}
Suppose that $\lambda \ntrianglerighteq \nu$.  
Then $x_\nu T_w y_{\lambda'}=0$ for all $w \in \mathscr{D}_\nu$. 
\end{lemma} 

\begin{proof}
The proof follows from Lemmas \ref{nicelemma2} and \ref{die2}.
\end{proof}

\begin{lemma}\label{die1}
Suppose $\lambda \ntrianglerighteq \nu$ and $\Theta:M^\lambda \rightarrow M^\nu$. Let $\hat{\Theta}$ denote the restriction of $\Theta$ to $S^\lambda$. Then $\hat{\Theta} =0$.
\end{lemma}

\begin{proof}
Recall that $S^\lambda$ is generated by $x_\lambda T_{w_\lambda} y_{\lambda'}$.  Then we may write
\[ \hat{\Theta}(x_\lambda T_{w_\lambda} y_{\lambda'}) = \sum_{w \in \mathscr{D}_\nu} f(w) x_{\nu} T_w y_{\lambda'}\]
for some $f(w) \in F$.  By Lemma \ref{die3}, $x_{\nu} T_w y_{\lambda'}=0$ for all $w \in \mathscr{D}_\nu$.
\end{proof}

The following theorem has been proved by Donkin~\cite{Donkin:tilting}, Proposition 10.4 and by Lyle and Mathas~\cite{LM:rowhoms}, Theorem 3.2; it will considerably simplify our later working.

\begin{theorem}\label{rowrem}
Suppose that $\lambda$ and $\mu$ are partitions of $n$ and that either $e\neq 2$ or $\lambda$ is 2-regular.  
\begin{itemize}
\item Suppose that $\lambda_1=\mu_1$.  Let $\bar{\lambda}=(\lambda_2,\lambda_3,\ldots)$ and $\bar{\mu}=(\mu_2,\mu_3,\ldots)$.  Then
\[\Hom_\h(S^\lambda,S^\mu) \cong_F \Hom_\h(S^{\bar{\lambda}},S^{\bar{\mu}}).\]
\item Suppose that $\lambda'_1=\mu'_1$.  Let $\bar{\lambda}=(\lambda_1-1,\lambda_2-1,\ldots)$ and $\bar{\mu}=(\mu_1-1,\mu_2-1,\ldots)$.  Then
\[\Hom_\h(S^\lambda,S^\mu) \cong_F \Hom_\h(S^{\bar{\lambda}},S^{\bar{\mu}}).\]
\end{itemize}
\end{theorem}

Finally in this section, we take a step backwards, from Specht modules to permutation modules.  Let $A \in \mathcal{T}(\lambda,\nu)$ be a row standard tableau and fix $d$ and $t$ with $d \geq 1$ and $0 \leq t <\mu_{d+1}$.  Let $\nu = \nu^{d,t}$.
We will consider the map $\psi_{d,t} \Theta_A:M^\lambda \rightarrow M^\nu$.  We write $\psi_{d,t} \Theta_A$ in terms of homomorphisms indexed by row standard $\lambda$-tableaux of type $\nu$.  To do so, we make use of the  Gaussian, or quantum, polynomials $\mat{\alpha}{\beta}$; a useful reference is~\cite{James:gl}.

\begin{definition}
Suppose $\alpha\geq 0$.  Let $[\alpha] \in F$ be defined by
\begin{align*} 
[0] & = 0, \\
[\alpha] &= 1 + q + q^2+\ldots+q^{\alpha-1}\text{ if } \alpha > 0 , \\
\intertext{and set}
[0]!&=1,\\ [\alpha]! &= [1][2]\ldots[\alpha] \text{ if } \alpha > 0.  \\
\intertext{If $\alpha \geq \beta \geq 0$, define}
\mat{\alpha}{\beta}&=\frac{[\alpha]!}{[\beta]![\alpha-\beta]!}.\\
\intertext{Then $\mat{\alpha}{\beta}$ can be shown to be a polynomial in $q$ with integer coefficients, and} 
\mat{\alpha}{\beta}&=\mat{\alpha-1}{\beta} +q^{\alpha-\beta}\mat{\alpha-1}{\beta-1}.\end{align*}
\end{definition}

\begin{lemma}\label{gausslem}
Fix $\alpha\geq \beta \geq 0$ and let $\mathcal{I}^\alpha_\beta=\{ I=(i_1,i_2,\ldots,i_\beta) \mid 1\leq i_i<i_2<\ldots<i_\beta\leq \alpha\}$.  For $I=(i_1,i_2,\ldots,i_\beta) \in \mathcal{I}^\alpha_\beta$, set \begin{align*}
G(I) &= \sum_{j=1}^\beta (\alpha-i_j -\beta+j) \intertext{and set}
\Sigma(\alpha,\beta) &= \sum_{I \in \mathcal{I}^\alpha_\beta} q^{G(I)}.\\  \intertext{Then}
\Sigma(\alpha,\beta) &= \mat{\alpha}{\beta}.\end{align*}
\end{lemma}

\begin{proof}
The proof is by induction on $\alpha$, the case $\alpha=0$ being trivial.  Suppose Lemma \ref{gausslem} holds for $\alpha-1$.  It is easy to see that
\begin{align*}
\Sigma(\alpha,\beta) &=\Sigma(\alpha-1,\beta) + q^{\alpha-\beta}\Sigma(\alpha-1,\beta-1) \\
&= \mat{\alpha-1}{\beta} + q^{\alpha-\beta}\mat{\alpha-1}{\beta-1}\\
\intertext{by the inductive hypothesis}
&= \mat{\alpha}{\beta}.
\end{align*}
\end{proof}

\begin{proposition}\label{comb}
Let $\lambda$ and $\mu$ be partitions of $n$ and choose $d$ and $t$ with $d \geq 1$ and $0 \leq t < \mu_{d+1}$.  Let $\nu = \nu^{d,t}$ and write $\bar{t}=\mu_{d+1}-t$.
Suppose $A \in \mathcal{T}(\lambda,\mu)$ is a row standard tableau.  Let $\mathcal{S} \subseteq \mathcal{T}(\lambda,\nu)$ be the set of row standard tableaux obtained by replacing $\bar{t}$ entries of $d+1$ in $A$ with $d$.  For $S \in \mathcal{S}$ and $i\geq 1$, suppose that $\beta_{i}$ entries were replaced in row $i$. Define $b_S \in F$ by
\[b_S= \prod_{i \geq 1} q^{x_i \beta_i} \mat{y_i}{\beta_i}\] where 
$x_i$ is the cardinality of the set $\{(k,j) \mid k>i \text{ and } A(k,j)=d\}$
and $y_i$ is the cardinality of the set $\{j \mid S(i,j)=d\}$.
Then \[ \psi_{d,t} \Theta_A = \sum_{S \in \mathcal{S}} b_S \Theta_S.\]
\end{proposition}
\begin{proof}  Let $R \in \mathcal{T}(\mu,\nu)$ be such that $R(d+1,b)=d$ for $b \leq \bar{t}$ and, for all other values of $(a,b)$, $R(a,b)=a$.  Note that the map $\psi_{d,t} \Theta_A=\Theta_R \Theta_A$ is completely determined by its action on $x_\lambda$.  
\[ \psi_{d,t} \Theta_A (x_\lambda) = x_{\nu} \Bigg( \sum_{R' \sim_r R} T_{1_{R'}}\Bigg) \Bigg(\sum_{A'\sim_r A}T_{1_{A'}}\Bigg).
\]
If $R'\sim_r R$ then $1_{R'} \in \sym{\mu}$. If $A' \sim_r A$ then $1_{A'} \in \mathscr{D}_\mu$ and hence $T_{1_{R'}}T_{1_{A'}} = T_{1_{R'}1_{A'}}$.

Choose $R'\sim_r R$ and $A'\sim_r A$.  Then $\mf{t}^\nu 1_{R'} 1_{A'}$ is formed by taking the tableau $\mf{t}^\mu 1_{A'}$ and raising $\bar{t}$ nodes from row $d+1$ of $\mf{t}^\mu 1_{A'}$ to the right end of row $d$. 
It is therefore row equivalant to a tableau $\mf{t}^\nu 1_{\bar{A}}$ where $\bar{A}$ is formed by replacing $\bar{t}$ of the entries of $d+1$ in $A'$ with $d$.  
Suppose that these entries were at nodes $(i,j_1),(i,j_2),\ldots,(i,j_{\bar{t}})$. 
Then the nodes moved were of value $\mf{t}^\lambda(i,j_1),\mf{t}^\lambda(i,j_2),\ldots,\mf{t}^\lambda(i,j_{\bar{t}})$.  For $1 \leq k \leq \bar{t}$, let $g(i,j_k)$ be equal to the cardinality of the the set $\{(x,y) \mid A'(x,y) =d \text{ and } \mf{t}^\lambda(x,y) > \mf{t}^\lambda(i,j_k)\}$.
The number of entries in row $d$ of $\mf{t}^\nu 1_{\bar{A}}$ which are greater than $\mf{t}^\lambda(i,j_k)$ is equal to $g(i,j_k)$. Set $G(\bar{A}) = \sum_{k=1}^{\bar{t}} g(i,j_k)$.  
Therefore 
\[x_\nu T_{1_{R'} 1_{A'}} = q^{G(\bar{A})} x_\nu T_{1_{\bar{A}}}.\]

Let $A' \sim_r A$.  Write $\bar{A} \rightarrow A'$ if $\bar{A} \in \mathcal{T}(\lambda,\nu)$ is formed by replacing $\bar{t}$ entries of $d+1$ in $A'$ by $d$.   Then
\[  x_{\nu}\sum_{R' \sim_r R} \sum_{A'\sim_r A} T_{1_{R'}1_{A'}} = x_\nu \sum_{A' \sim_r A} \sum_{\bar{A} \rightarrow A'} q^{G(\bar{A})}T_{1_{\bar{A}}}.\]
The result then follows from Lemma \ref{gausslem}.
\end{proof}

\section{Trivial submodules of Specht modules}
Let $\h = \h_{F,q}(\sym{n})$ where $F$ is a field of characteristic $p \geq 0$ and define $e > 1$ to be minimal such that $1+q+\ldots+q^{e-1}=0$.  Since $\h$ is semisimple if $e = \infty$, we may assume that $e$ is finite; our results trivially hold if $e = \infty$. 

In this short section, we determine which Specht modules contain a submodule isomorphic to the trivial module $S^{(n)}$ by calculating the homomorphism spaces $\Hom_\h(S^{(n)},S^{\mu})$ for all $\mu$.  This generalises the result of James~\cite{James}, Theorem 24.4 for the symmetric groups; our approach is an exact analogue.  Naturally, $\dim(\Hom_\h(S^{(n)},S^{\mu})) \leq 1$ for all $\mu$ and  $\dim(\Hom_\h(S^{(n)},S^{\mu}))= 1$ if and only if $S^\mu$ has a submodule isomorphic to $S^{(n)}$.

\begin{lemma}
Suppose that $p=0$.  
Let $\alpha\geq 0,\beta \geq 1$.  Then 
\[\mat{\alpha+1}{1},\mat{\alpha+2}{2},\ldots,\mat{\alpha+\beta}{\beta} \]
are all zero in $F$ if and only $e \mid \alpha + 1$ and $\beta < e$.
\end{lemma}

\begin{proof}
We have that $\mat{\alpha+1}{1}=0$ if and only if $e \mid \alpha +1$; and if $e \mid \alpha+1$ then $\mat{\alpha+\gamma}{\gamma}=0$ for all $\gamma <e$.  Suppose $e \mid \alpha+1$ and consider $\mat{\alpha+e}{e}$. Clearly it is zero if and only if $\frac{[\alpha+1]}{[e]} =0$.  But 
\[\frac{[\alpha+1]}{[e]} = 1 + q^e + q^{2e}+\ldots q^{\alpha+1 -e} = \frac{\alpha+1}{e} \neq 0.\]
\end{proof}

\begin{lemma}[{\cite{James:gl}, Theorem 19.5}] \label{gausszero}
Suppose $p > 0$.
For each non--negative integer $b$, write $b = b^\ast e + b'$ where $0 \leq b' <e$, and define $\ell_p(b)$ to be minimal such that $b < p^{\ell_p(b)}$.
Let $\alpha\geq 0,\beta \geq 1$.  Then 
\[\mat{\alpha+1}{1},\mat{\alpha+2}{2},\ldots,\mat{\alpha+\beta}{\beta} \]
are all zero in $F$ if and only 
\[\alpha  \equiv -1 \mod ep^{\ell_p(\beta^\ast)}.\]
\end{lemma}

\begin{theorem}\label{maphappens}
Take $\mu = (\mu_1,\mu_2,\ldots,\mu_l)$ to be a partition of $n$ with exactly $l$ parts.

Suppose $p=0$.
The Specht module $S^\mu$ has a submodule isomorphic to the trivial $\h$-module  $S^{(n)}$ if and only if $\mu=(n)$ or $\mu=(\mu_1,(e-1)^{l-2},\mu_l)$ where $e \mid \mu_1+1$. 

Suppose $p>0$
The Specht module $S^\mu$ has a submodule isomorphic to the trivial $\h$-module  $S^{(n)}$ if and only if for $1 \leq i <l$, $\mu_i  \equiv -1 \mod ep^{z_i}$ where $z_i = \ell_p((\mu_{i+1})^\ast)$.

\end{theorem}

\begin{proof}
Suppose $\hat{\Theta}:S^{(n)} \rightarrow M^{\mu}$.  Then $\hat{\Theta}$ is a linear multiple of the map $\hat{\Theta}_A$ where $A$ is the unique semistandard $(n)$-tableau of type $\mu$.  Recall that $\Image{\hat{\Theta}}\subseteq S^\mu$ if and only if $\psi_{d,t} \hat{\Theta} = 0$ for $1 \leq d <l$ and $0 \leq t <\mu_{d+1}$.  Fix $d$ with $1 \leq d <l$.  For $0 \leq t <\mu_{d+1}$, let $S \in \mathcal{T}_0(\lambda,\nu^{d,t})$ be the tableau obtained by replacing the first $\mu_{d+1}-t$ entries of $d+1$ in $A$ by $d$. 
By Proposition \ref{comb},
\[ \psi_{d,t}\hat{\Theta}= \mat{\mu_d+\mu_{d+1}-t}{\mu_{d+1}-t} \hat{\Theta}_S.\] 
Then $\psi_{d,t}\hat{\Theta}= 0$ for all $0 \leq t <\mu_{d+1}$ if and only if 
$\mat{\mu_d + \beta}{\beta} = 0$ for all $1 \leq \beta \leq \mu_{d+1}$.
\end{proof}

The following theorems are $q$-analogues of the Carter--Payne theorem, where we move nodes between adjacent rows.

\begin{theorem}\label{CProw1}
Suppose that $\lambda$ and $\mu$ are partitions of $n$ such that 
\[ \lambda_i = \begin{cases} \mu_i +\gamma & \text{if } i = a, \\
\mu_i -\gamma & \text{if } i = a+1, \\
\mu_i &\text{otherwise},
\end{cases}\]
for some positive integers $a$ and $\gamma$, and that $\lambda$ is 2-regular if $e=2$. If $p=0$ \begin{align*}
\dim(\Hom_\h(S^\lambda,S^\mu)) &= \begin{cases} 1 & \text{if } \mu_d -\mu_{d+1} + \gamma \equiv -1 \mod e \text{ and } \gamma < e, \\
0 & \text{otherwise}. \end{cases}\\
\intertext{If $p >0$,}
\dim(\Hom_\h(S^\lambda,S^\mu)) &= \begin{cases} 1 & \text{if } \mu_d -\mu_{d+1} + \gamma \equiv -1 \mod e p^{\ell_p(\gamma^\ast)}, \\
0 & \text{otherwise}. \end{cases} \end{align*}
\end{theorem}

\begin{proof}
The proof follows from Theorems \ref{maphappens} and Theorem \ref{rowrem}.
\end{proof}

When $e=2$ and $\lambda$ is 2-regular, we must be a little more circumspect.  The following result can easily be deduced from the proof of~\cite{LM:rowhoms}, Theorem 3.1.

\begin{theorem}\label{rowrem2}
Suppose that $\lambda$ and $\mu$ are partitions of $n$. 

Suppose that $\lambda_1=\mu_1$.  Let $\eta=(\lambda_2,\lambda_3,\ldots)$ and $\xi=(\mu_2,\mu_3,\ldots)$.  
For $A \in \mathcal{T}_0(\lambda,\mu)$, define $\bar{A} \in \mathcal{T}_0(\eta,\xi)$ by setting $\bar{A}(i,j) = A(i+1,j)$, and note that this gives a bijection between  $\mathcal{T}_0(\lambda,\mu)$ and $\mathcal{T}_0(\eta,\xi)$.  Suppose $\hat{\Theta}:S^\lambda \rightarrow M^\mu$ and $\hat{\Theta}':S^{\eta} \rightarrow M^{\xi}$ are such that \[ \hat{\Theta} = \sum_{A \in \mathcal{T}_0(\lambda,\mu)} f(A) \hat{\Theta}_A \qquad \qquad \hat{\Theta}' = \sum_{A \in \mathcal{T}_0(\lambda,\mu)} f(A) \hat{\Theta}_{\bar{A}} \] for some $f(A) \in F$.
Then $\Image(\hat{\Theta}) \subseteq S^{\mu}$ if and only if $\Image(\hat{\Theta}') \subseteq S^{\xi}$.

There is a similar theorem concerning column removal.
\end{theorem}

\begin{corollary}
Suppose that $\lambda$ and $\mu$ are partitions of $n$ such that 
\[ \lambda_i = \begin{cases} \mu_i +\gamma & \text{if } i = a, \\
\mu_i -\gamma & \text{if } i = a+1, \\
\mu_i &\text{otherwise},
\end{cases}\]
for some positive integers $a$ and $\gamma$.  
If $p=0$, \begin{align*} 
\dim(\Hom_\h(S^\lambda,S^\mu)) &\geq 1  \text{ if } \mu_d -\mu_{d+1} + \gamma \equiv -1 \mod e \text{ and } \gamma < e. \\
\intertext{If $p >0$,}
\dim(\Hom_\h(S^\lambda,S^\mu)) &\geq 1  \text{ if } \mu_d -\mu_{d+1} + \gamma \equiv -1 \mod e p^{\ell_p(\gamma^\ast)}. \end{align*}
\end{corollary}

\numberwithin{equation}{subsection}
\section{One node Carter--Payne homomorphisms}
\subsection{Backround}
We now concentrate on pairs of partitions $\lambda$ and $\mu$, where $\lambda$ is formed from $\mu$ by raising one node. 
By Theorem \ref{rowrem} and Theorem \ref{rowrem2}, the following two theorems are equivalent.
\begin{theorem}\label{mainextra}
Suppose that
\begin{align*}
\mu&=(\mu_1,\ldots,\mu_{a-1},\mu_a,\mu_{a+1},\ldots,\mu_{b-1},\mu_b,\mu_{b+1},\ldots,\mu_r), \\
\lambda&=(\mu_1,\ldots,\mu_{a-1},\mu_a+1,\mu_{a+1},\ldots,\mu_{b-1},\mu_b-1,\mu_{b+1},\ldots,\mu_r),
\end{align*}
are partitions of $n$.  If $e \mid \mu_a - \mu_b +b-a+1$ then there exists $0 \neq \hat{\Theta}:S^\lambda \rightarrow S^\mu$, where
\[\hat{\Theta} = \sum_{A\in \mathcal{T}_0(\lambda,\mu)}f(A) \hat{\Theta}_A \]
for some $f(A) \in F$.  
\end{theorem}

\begin{theorem}\label{main}
Suppose that
\begin{align*}
\mu&=(\mu_1,\mu_2,\ldots,\mu_s,1), \\
\lambda&=(\mu_1+1,\mu_2,\ldots,\mu_s),
\end{align*}
are partitions of $n$.  If $e \mid \mu_1+s$ then there exists $0 \neq \hat{\Theta}:S^\lambda \rightarrow S^\mu$, where
\[\hat{\Theta} = \sum_{A\in \mathcal{T}_0(\lambda,\mu)}f(A) \hat{\Theta}_A \]
for some $f(A) \in F$.    
\end{theorem}

We shall give a constructive proof of Theorem \ref{main}. A direct proof of Theorem \ref{mainextra} would be very similar; we choose to consider Theorem \ref{main} mainly for convenience of notation.
Henceforth in this section, we fix partitions of $n$, 
\begin{align*}
\mu&=(\mu_1,\mu_2,\ldots,\mu_s,1), \\
\lambda&=(\mu_1+1,\mu_2,\ldots,\mu_s).
\end{align*}

Consider $\{\mathcal{T}_0(\lambda,\mu) \}$.  The tableaux $A$ in this set are determined by the following properties.  For $1 \leq a \leq s$,
\begin{itemize}
\item $A(a,b)=a$ for $b \neq \lambda_a$.
\item Write $A(a,\lambda_a)=i_a$.  Then $\{i_1,i_2,\ldots,i_s\} = \{2,3,\ldots,s+1\}$ where $i_a \geq a$, and if $\lambda_a=\lambda_{a+1}$ then $i_a<i_{a+1}$.
\end{itemize}
Hence for $A\in \mathcal{T}_0(\lambda,\mu)$, we will write $A=(\mu:i_1,i_2,\ldots,i_s)$.  For $2 \leq a' \leq s+1$, define $r(a')$ by $i_{r(a')}=a'$.

We now fix a map $\hat{\Theta}:S^\lambda \rightarrow M^\mu$, setting
\[\hat{\Theta} = \sum_{A\in \mathcal{T}_0(\lambda,\mu)}f(A) \hat{\Theta}_A \]
for some $f(A) \in F$.  We will write $f(A)=f(\mu:i_1,i_2,\ldots,i_s)$.

\begin{lemma}\label{onlyone}
For $1\leq d \leq s$, write $\psi_d = \psi_{d,\mu_{d+1}-1}:M^\mu \rightarrow M^{\nu^{d,\mu_{d+1}-1}}$.  Then $\Image(\hat{\Theta}) \subseteq S^\mu$ if and only if $\psi_{d} \hat{\Theta} = 0$ for all $1\leq d\leq s$.
\end{lemma}

\begin{proof}
For $1\leq d \leq s$, choose $t$ with $0 \leq t < \mu_{d+1}-1$.  Then $\lambda \ntrianglerighteq \nu^{d,t}$ so that by Lemma \ref{die1}, $\psi_{d,t}\hat{\Theta} =0$.   The result then follows from Corollary \ref{cor:compose}.
\end{proof}

Our aim is therefore to rewrite each $\psi_{d} \hat{\Theta}$ in terms of semistandard homomorphisms, and to deduce necessary and sufficient conditions for the coefficient of each semistandard homomorphism to be zero.  We begin with some preliminary results before discussing the maps $\hat{\Theta}_S$, where $S \in \mathcal{T}(\lambda,\nu^{d,t})$.  We are then able to rewrite the maps $\psi_{d} \hat{\Theta}_A$ in terms of semistandard homomorphisms.  Finally, in Proposition \ref{sumconds}, we describe three straightforward conditions that specify when $\psi_{d} \hat{\Theta} = 0$ for all $1\leq d\leq s$.

We remark that while the results concerning the manipulation of the maps $\hat{\Theta}_S$ tend to be reasonably simple, the only proofs that we have been able to discover have usually been somewhat involved; in particular, they are a lot more complicated than the corresponding proofs for the symmetric groups.

\subsection{Preliminary results}
For $1 \leq d \leq s$, we define compositions $\nu=\nu(d)$, and for $1 \leq d < s$ we define compositions $\sigma=\sigma(d)$, as follows.
\[\nu_i = \begin{cases}
\mu_i +1 & \text{if } i=d,\\ 
\mu_i-1 & \text{if } i=d+1,\\
\mu_i & \text{otherwise,} \end{cases} \qquad
\sigma_i = \begin{cases}
\lambda_i +1 & \text{if }i=d,\\ 
\lambda_i-1 & \text{if }i=d+1,\\
\lambda_i & \text{otherwise}. \end{cases} \]
For convenience, we introduce two more items of notation.  For $1 \leq x \leq y \leq n$, write
\begin{align*}
\up{x}{y} &= (I + T_x +T_xT_{x+1} + \ldots + T_x T_{x+1} \ldots T_{y-1}), \\
\down{x}{y} & = (I + T_{y-1} + T_{y-1}T_{y-2}+\ldots+ T_{y-1}T_{y-2}\ldots T_x),
\end{align*}
where $I$ denotes the identity element of $\h$.

\begin{lemma}\label{nicelemma1}
For $1 \leq d < s$,
\[x_\sigma \up{\mu_1+\ldots+\mu_d+2}{\mu_1+\ldots+\mu_{d+1}+1}T_{w_\lambda}y_{\lambda'} =0. \] 
\end{lemma}
\begin{proof}
Let $R \in \mathcal{T}(\lambda,\sigma)$ be defined by
\[ R(a,b) = \begin{cases} d & \text{if } (a,b) = (d+1,1), \\
a & \text{otherwise}. \end{cases}\]
From the definition of $S^\lambda$ in Theorem \ref{altdefn},
\[0 = \hat{\Theta}_R(x_\lambda T_{w_\lambda} y_{\lambda'}) = x_\sigma \up{\mu_1+\ldots+\mu_d+2}{\mu_1+\ldots+\mu_{d+1}+1}T_{w_\lambda}y_{\lambda'}.\]
\end{proof}

\begin{lemma}\label{separation}
Let $S \in \mathcal{T}(\lambda,\nu)$ and suppose there exists $i$ such that $S(a,b)  \leq S(a',b')$ whenever $a \leq i$ and $a' > i$.  Define tableaux $S^t$ and $S^b$ as follows
\begin{align*}
S^t(a,b) &= \begin{cases} S(a,b) & \text{if } a \leq i, \\
a & \text{if } a > i, \end{cases} &
  S^b(a,b) &= \begin{cases} a & \text{if } a \leq i, \\
S(a,b) & \text{if } a > i. \end{cases}\\
\intertext{and set} 
h^t & = \sum_{S' \sim_r S^t} T_{1_{S'}} &
h^b & = \sum_{S' \sim_r S^b} T_{1_{S'}}.
\end{align*}Then
\[  \Bigg(\sum_{S'\sim_r S} T_{1_{S'}}\Bigg) = h^t h^b.  \]
Note that $h^t$ lies within the subalgebra generated by $\{T_k \mid k \leq \lambda_1+\ldots+\lambda_i-1 \}$
and $h^b$ lies within the subalgebra of $\h$ generated by $\{T_k \mid k \geq \lambda_1+\ldots+\lambda_i+1 \}$, so that $h^t$ and $h^b$ commute.  Note also that while $S^t$ and $S^b$ are of shape $\lambda$, we do not need to specify their type.
\end{lemma}

\begin{proof}
The lemma follows from the definition of the permutations $1_{S'}$.
\end{proof}

The following lemma may be proved by induction; we leave the proof as an exercise for the reader.

\begin{lemma}\label{mess}
Choose $1\leq x \leq y \leq n$.  Then
\begin{align*}
T_{y-1}& T_{y-2}\ldots T_x \up{x}{y}\\
= &\, q^{y-x} I + q^{y-x} T_{y-1} +q^{y-x-1}T_{y-1} T_{y-2} + q^{y-x-2} T_{y-1} T_{y-2} T_{y-3} + \ldots +q T_{y-1} T_{y-2} \ldots T_x \\
& +q^{y-x-2}(q-1)T_{y-2}T_{y-1} T_{y-2} \\
& +q^{y-x-3}(q-1)(T_{y-3}T_{y-1} T_{y-2} T_{y-3} + T_{y-3}T_{y-2}T_{y-1} T_{y-2} T_{y-3})\\
& \ldots\\
& +(q-1)(T_xT_{y-1}T_{y-2}\ldots T_x + T_x T_{x+1} T_{y-1}T_{y-2}\ldots T_x +\ldots +T_x\ldots T_{y-2}T_{y-1} T_{y-2}\ldots T_x).
\end{align*}
Therefore, if we choose $d$ with $1 \leq d <s$ and take 
\begin{align*}
x&= \mu_1 + \ldots+\mu_d +2, \\
y&= \mu_1+\ldots+\mu_d+\mu_{d+1} +1, \\
y \leq z & \leq n,
\end{align*}
then \[x_\nu T_{z-1}T_{z-2}\ldots T_{x}\up{x}{y} = q^{y-x}x_\nu T_{z-1}T_{z-2}\ldots T_{y} \down{x}{y}.\]
\end{lemma}
This completes our preliminary results.

\subsection{Manipulation of maps}
Before studying the maps $\psi_d\hat{\Theta}_A:S^\lambda \rightarrow M^\nu$, we collect together some information about the maps $\hat{\Theta}_S:S^\lambda \rightarrow M^\nu$, where $S \in \mathcal{T}(\lambda,\nu)$. 

Choose $1\leq d \leq s$ and consider $\{\mathcal{T}_0(\lambda,\nu) \}$.  The tableaux $S$ in this set are determined by the following properties.  For $1 \leq a \leq s$,
\begin{itemize}
\item $S(a,b)=a$ for $b \neq \lambda_a$.
\item Write $(a,\lambda_a)=j_a$.  Then $\{j_1,j_2,\ldots,j_{d-1}\} = \{2,3,\ldots,d\},\{j_{d+1},j_{d+2},\ldots,j_s\}=\{d+2,d+3,\ldots,s+1\}$ and $j_d=d$, where $j_a \geq a$ and if $\lambda_a=\lambda_{a+1}$ then $j_a<j_{a+1}$.
\end{itemize}
If $S \in \mathcal{T}(\lambda,\nu)$ satisfies all of the conditions above, except possibly the condition that $j_a<j_{a+1}$ whenever $\lambda_a=\lambda_{a+1}$, we will write $S=(\nu:j_1,j_2,\ldots,j_s)$. 
For $1 \leq d \leq s$, define $\check{r}(d)$ by specifying that $i_{\check{r}(d)}=d$ and $\check{r}(d)<d$.

\begin{lemma}\label{swap1}
Let $A=(\mu:i_1,i_2,\ldots,i_s) \in \mathcal{T}_0(\lambda,\mu)$ be such that $i_d=d+1$ and $A(d+1,1)=d+1$.  Let $S \in \mathcal{T}(\lambda,\nu)$ be the row standard tableau formed by replacing the entry $A(d+1,1)$ with $d$, and let $U \in\mathcal{T}(\lambda,\nu)$ be the row standard tableau formed by replacing the entry $A(d,\lambda_d)$ with $d$.  Then
\[ \hat{\Theta}_S = - [\mu_{d+1}-1] \hat{\Theta}_{U}.\] 
\end{lemma}
\begin{proof}
Let $R \in \mathcal{T}(\sigma,\nu)$ be the tableau formed by setting
\[ R(a,b) = \begin{cases} d+1 & \text{if } (a,b) = (d,\lambda_d+1), \\
U(a,b+1) & \text{if } a=d+1, \\
U(a,b) & \text{otherwise}. 
\end{cases}\]
Then \begin{multline*}
 \Bigg(\sum_{R'\sim_r R} T_{1_{R'}}\Bigg)\up{\mu_1+\ldots+\mu_d+2}{\mu_1+\ldots+\mu_{d+1}+1} = 
\sum_{S'\sim_r S} T_{1_{S'}} \\
+\up{\mu_1+\ldots+\mu_d+2}{\mu_1+\ldots+\mu_{d+1}} \Bigg(\sum_{U' \sim U} T_{1_{U'}}\Bigg).
\end{multline*}
Note that the map $\hat{\Theta}_S$ is completely determined by its action on $x_\lambda T_{w_\lambda} y_{\lambda'}$.  
\begin{align*} 
\hat{\Theta}_{S} (x_\lambda T_{w_\lambda} y_{\lambda'}) &= x_\nu \Bigg(\sum_{S'\sim_r S} T_{1_{S'}}\Bigg)  T_{w_\lambda} y_{\lambda'}\\
&= x_\nu \Bigg(\sum_{R'\sim_r R} T_{1_{R'}}\Bigg)\up{\mu_1+\ldots+\mu_d+2}{\mu_1+\ldots+\mu_{d+1}+1} T_{w_\lambda} y_{\lambda'} \\
&\qquad - x_\nu \up{\mu_1+\ldots+\mu_d+2}{\mu_1+\ldots+\mu_{d+1}} \Bigg(\sum_{U'\sim_r U} T_{1_{U'}}\Bigg) T_{w_\lambda} y_{\lambda'} \\
& = h' x_\sigma \up{\mu_1+\ldots+\mu_d+2}{\mu_1+\ldots+\mu_{d+1}+1} T_{w_\lambda} y_{\lambda'} \\
&\qquad - x_\nu \up{\mu_1+\ldots+\mu_d+2}{\mu_1+\ldots+\mu_{d+1}} \Bigg(\sum_{U'\sim_r U} T_{1_{U'}}\Bigg) T_{w_\lambda} y_{\lambda'} \\
\intertext{for some $h' \in \h$}
&= - [\mu_{d+1}-1]x_\nu \Bigg(\sum_{U'\sim_r U} T_{1_{U'}}\Bigg) T_{w_\lambda} y_{\lambda'} \\
\intertext{by Lemma \ref{nicelemma1}, and noting that $\up{\mu_1+\ldots+\mu_d+2}{\mu_1+\ldots+\mu_{d+1}} \subseteq \sym{\nu}$}
&= - [\mu_{d+1}-1] \hat{\Theta}_{U}(x_\lambda T_{w_\lambda} y_{\lambda'}).
\end{align*}
\end{proof}

\begin{lemma}\label{kill2}
Suppose that $A=(\mu:i_1,i_2,\ldots,i_s) \in \mathcal{T}_0(\lambda,\mu)$ is such that $i_d=d+1$. Let $U \in\mathcal{T}(\lambda,\nu)$ be the row standard tableau formed by replacing the entry $A(d,\lambda_d)$ with $d$.  Unless $\mu_{d-1}=\mu_d$ and $i_{d-1}=d$, the tableau $U$ is semistandard.  If $\mu_{d-1}=\mu_d$ and $i_{d-1}=d$ then  $\hat{\Theta}_{U}=0$.
\end{lemma}
\begin{proof}
The first part of the lemma is obvious.  Suppose that $\mu_{d-1}=\mu_d$ and $i_{d-1}=d$; note that $i_a \leq d-1$ for $a \leq d-2$.  Choose $U' \sim_r U$. Using the same technique as Lemma \ref{separation}, it is possible to write
\[1_{U'} = v^t v^b v\]
where 
\begin{align*} v^t &\text{ lies in the subgroup } \sym{(1,2,\ldots,\mu_1+\ldots+\mu_{d-2}+1)}, \\
v^b &\text{ lies in the subgroup }\sym{(\mu_1+\ldots+\mu_d+2,\ldots,s)},\\
v &\text{ lies in the subgroup }\sym{(\mu_1+\ldots+\mu_{d-2}+2,\ldots,\mu_1+\ldots+\mu_{d-1}+1)}.
\end{align*}

Then $v \in \sym{\lambda}$ so $T_v T_{w_\lambda} = T_{v w_\lambda}$, and from the properties of $U$, $\mf{t}^\nu v w_\lambda$ is a tableau such that row $d$ contains two numbers from the same column of $\mf{t}_\lambda$.  Hence $x_\nu T_{1_{U'}} T_{w_\lambda}$ is equal to a sum 
\[x_\nu \sum_{w \in\mathscr{D}_\nu} f(w) T_{w}\]
for some $f(w) \in F$ where each tableau $\mf{t}^\nu w$ has the property that row $d$ contains two entries which come from the same column of $\mf{t}_\lambda$.  It therefore follows from Lemma \ref{nicelemma2} that $x_\nu T_{w} y_{\lambda'} =0$.
\end{proof}

\begin{lemma}\label{kill1}
Let $A=(\mu:i_1,i_2,\ldots,i_s) \in \mathcal{T}_0(\lambda,\mu)$ be such that $i_d=d$ and $A(d+1,1)=d+1$.  Let $S \in \mathcal{T}(\lambda,\nu)$ be the row standard tableau formed by replacing the entry $A(d+1,1)$ with $d$.  Then $\hat{\Theta}_S=0$.
\end{lemma}

\begin{proof}
Let $R \in \mathcal{T}(\sigma,\nu)$ be the tableau formed by setting
\[ R(a,b) = \begin{cases} d & \text{if } (a,b) = (d,\lambda_d+1), \\
S(a,b+1) & \text{if } a=d+1, \\
S(a,b) & \text{otherwise}. 
\end{cases}\]
Then \begin{align*}
\sum_{S'\sim_r S} T_{1_{S'}} 
&= \Bigg(\sum_{R'\sim_r R} T_{1_{R'}}\Bigg)\up{\mu_1+\ldots+\mu_d+2}{\mu_1+\ldots+\mu_{d+1}+1}\\
\intertext{and therefore}
\hat{\Theta}_{S} (x_\lambda T_{w_\lambda} y_{\lambda'}) &= x_\nu \Bigg(\sum_{S'\sim_r S} T_{1_{S'}}\Bigg)  T_{w_\lambda} y_{\lambda'}\\
&= x_\nu  \Bigg(\sum_{R'\sim_r R} T_{1_{R'}}\Bigg)\up{\mu_1+\ldots+\mu_d+2}{\mu_1+\ldots+\mu_{d+1}+1} T_{w_\lambda} y_{\lambda'}\\
&= h' x_\sigma \up{\mu_1+\ldots+\mu_d+2}{\mu_1+\ldots+\mu_{d+1}+1} T_{w_\lambda} y_{\lambda'}\\
\intertext{for some $h'\in \h$}
&=0 \end{align*}
by Lemma \ref{nicelemma1}.
\end{proof}

\begin{lemma}\label{swap2}
Let $A=(\mu:i_1,i_2,\ldots,i_s) \in \mathcal{T}_0(\lambda,\mu)$ be such that $i_{d+1}=d+1$ and $i_d \neq d$.  Let $S \in \mathcal{T}(\lambda,\nu)$ be the row standard tableau formed by replacing the entry $A(d+1,1)$ with $d$. 
Let $U = (\nu:i_1,i_2,\ldots,i_{d-1},d,i_{d},i_{d+2}\ldots,i_s) \in \mathcal{T}(\lambda,\nu)$.
Then $\hat{\Theta}_S=-q^{\mu_{d+1}-1}\hat{\Theta}_{U}$.
\end{lemma}
\begin{proof} Note that $i_j = j$ for $d+2 \leq j <i_d$.
Define $\lambda$-tableaux $W$ and $\bar{S}$ by 
\[ W(a,b) = \begin{cases} S(a,b) & \text{if }a <d \text{ or } a \geq i_d, \\
a & \text{otherwise}, \end{cases} 
\qquad \bar{S}(a,b) = \begin{cases} a & \text{if }a <d \text{ or } a \geq i_d, \\
S(a,b) & \text{otherwise}, \end{cases}\]
 and set 
\[ \bar{h} =  \sum_{W'\sim_r W} T_{1_{W'}}.\]
By Lemma \ref{separation},
\[ \sum_{S'\sim_r S} T_{1_{S'}} =  \Bigg( \sum_{S'\sim_r \bar{S}}  T_{1_{S'}} \Bigg) \bar{h}\]where these two terms commute.
Furthermore,
\begin{multline*} \Bigg( \sum_{S'\sim_r \bar{S}} T_{1_{S'}}\Bigg)
= T_{\mu_1+\ldots+\mu_{i_d-1}}T_{\mu_1+\ldots+\mu_{i_d-1}-1} \ldots T_{\mu_1+\ldots+\mu_{d-1}+2} \\
\up{\mu_1+\ldots+\mu_d+2}{\mu_1+\ldots+\mu_{d+1}+1} \down{\mu_1+\ldots+\mu_{d-1}+2}{\mu_1+\ldots+\mu_{d}+1}.\end{multline*}
Let $R \in \mathcal{T}(\sigma,\nu)$ be the tableau formed by setting
\[ R(a,b) = \begin{cases} i_d & \text{if } (a,b) = (d,\lambda_d+1), \\
S(a,b+1) & \text{if } a=d+1, \\
S(a,b) & \text{otherwise}. 
\end{cases}\]
Then \begin{multline*}
\Bigg(\sum_{R'\sim_r R} T_{1_{R'}}\Bigg)\up{\mu_1+\ldots+\mu_d+2}{\mu_1+\ldots+\mu_{d+1}+1} = 
\sum_{S'\sim_r S} T_{1_{S'}} \\
+ T_{\mu_1+\ldots+\mu_{i_d-1}}T_{\mu_1+\ldots+\mu_{i_d-1}-1} \ldots T_{\mu_1+\ldots+\mu_{d-1}+2} \up{\mu_1+\ldots+\mu_{d-1}+2}{\mu_1+\ldots+\mu_d+1} \bar{h}.
\end{multline*}
Therefore
\begin{align*} 
\hat{\Theta}_{S} (x_\lambda T_{w_\lambda} y_{\lambda'}) &= x_\nu \Bigg(\sum_{S'\sim_r S} T_{1_{S'}}\Bigg)  T_{w_\lambda} y_{\lambda'}\\
&= x_\nu \Bigg(\sum_{R'\sim_r R} T_{1_{R'}}\Bigg)\up{\mu_1+\ldots+\mu_d+2}{\mu_1+\ldots+\mu_{d+1}+1} T_{w_\lambda} y_{\lambda'} \\
&\qquad - x_\nu  T_{\mu_1+\ldots+\mu_{i_d-1}}T_{\mu_1+\ldots+\mu_{i_d-1}-1} \ldots T_{\mu_1+\ldots+\mu_{d-1}+2} \\
&\qquad \qquad \up{\mu_1+\ldots+\mu_{d-1}+2}{\mu_1+\ldots+\mu_d+1} \bar{h} T_{w_\lambda} y_{\lambda'}\\
& = h' x_\sigma \up{\mu_1+\ldots+\mu_d+2}{\mu_1+\ldots+\mu_{d+1}+1} T_{w_\lambda} y_{\lambda'} \\
&\qquad - x_\nu  T_{\mu_1+\ldots+\mu_{i_d-1}}T_{\mu_1+\ldots+\mu_{i_d-1}-1} \ldots T_{\mu_1+\ldots+\mu_{d-1}+2} \\
&\qquad \qquad \up{\mu_1+\ldots+\mu_{d-1}+2}{\mu_1+\ldots+\mu_d+1} \bar{h} T_{w_\lambda} y_{\lambda'} \\
\intertext{for some $h' \in \h$}
&=-  x_\nu  T_{\mu_1+\ldots+\mu_{i_d-1}}T_{\mu_1+\ldots+\mu_{i_d-1}-1} \ldots T_{\mu_1+\ldots+\mu_{d-1}+2} \\
&\qquad \up{\mu_1+\ldots+\mu_{d-1}+2}{\mu_1+\ldots+\mu_d+1} \bar{h} T_{w_\lambda} y_{\lambda'}\\
\intertext{by Lemma \ref{nicelemma1}}
&=-  q^{\mu_{d+1}-1}x_\nu  T_{\mu_1+\ldots+\mu_{i_d-1}}T_{\mu_1+\ldots+\mu_{i_d-1}-1} \ldots T_{\mu_1+\ldots+\mu_{d}+1} \\
&\qquad \down{\mu_1+\ldots+\mu_{d-1}+2}{\mu_1+\ldots+\mu_d+1} \bar{h} T_{w_\lambda} y_{\lambda'}\\
\end{align*}
by Lemma \ref{mess}.
But it is straightforward to see that if we define $U$ as in Lemma \ref{swap2} then 
\begin{multline*}  \sum_{U'\sim_r U} T_{1_{U'}} = T_{\mu_1+\ldots+\mu_{i_d-1}}T_{\mu_1+\ldots+\mu_{i_d-1}-1} \ldots T_{\mu_1+\ldots+\mu_{d}+1} \\
\down{\mu_1+\ldots+\mu_{d-1}+2}{\mu_1+\ldots+\mu_d+1} \bar{h} \end{multline*}
and so $\hat{\Theta}_S = -q^{\mu_{d+1}-1} \hat{\Theta}_{U}$.
\end{proof}

\begin{lemma} \label{littleswap}
Suppose $c,d$ are such that $1 \leq d < c \leq s$ and $\lambda_c = \lambda_{c+1}$.  Suppose $U=(\nu: i_1,\ldots,i_s),V= (\nu: i'_1,\ldots,i'_s) \in \mathcal{T}(\lambda,\nu)$ are row standard tableaux with the following properties.
\begin{align*}
i_j &\leq c \text{ for }j <c, & &\\
i_c & =k \text{ for some }k>c+1, & i'_j &=i_j \text{ for } j \neq c,c+1,\\
i_j& = j \text{ for } c+1 \leq j <k, & i'_c & =c+1, \\
i_j &\geq k \text{ for }j \geq k, & i'_{c+1}& = k. 
\end{align*}
Then \[ \hat{\Theta}_U = - \hat{\Theta}_V.\]
\end{lemma}

\begin{proof}
Define the $\lambda$-tableau $W$ by 
\[ W(a,b) = \begin{cases} U(a,b) & \text{if }a <c \text{ or } a \geq k, \\
a & \text{otherwise}, \end{cases} \] and set 
\[ \bar{h} =  \sum_{W'\sim_r W} T_{1_{W'}}.\]
Define the $\lambda$-tableaux $\bar{U},\bar{V}$ by 
\[ \bar{U}(a,b) = \begin{cases} a & \text{if } a \neq c,c+1, \\
 a & \text{if } a = c,c+1 \text{ and } b < \mu_c, \\
 c+2 & \text{if } a = c \text{ and } b = \mu_c, \\
c +1 & \text{if } a = c+1 \text{ and } b = \mu_c, 
\end{cases} \qquad 
 \bar{V}(a,b) = \begin{cases} a & \text{if } a \neq c,c+1, \\
 a & \text{if } a = c,c+1 \text{ and } b < \mu_c, \\
 c+1 & \text{if } a = c \text{ and } b = \mu_c, \\
c +2 & \text{if } a = c+1 \text{ and } b = \mu_c. 
\end{cases} \]
Let $w \in \sym{n}$ be the permutation $(\mu_1+\ldots+\mu_{c+1}+1, \mu_1+\ldots+\mu_{c+1}+2,\ldots,\mu_1+\ldots+\mu_{k-1}-1)$.  Then
\begin{align*}
\hat{\Theta}_{U}(x_\lambda T_{w_\lambda} y_{\lambda'}) &= x_\nu \bar{h}T_w \Bigg( \sum_{U'\sim_r \bar{U}} T_{1_{U'}} \Bigg) T_{w_\lambda} y_{\lambda'}, \\
\hat{\Theta}_{V}(x_\lambda T_{w_\lambda} y_{\lambda'}) &= x_\nu \bar{h}T_w \Bigg( \sum_{V'\sim_r \bar{V}} T_{1_{V'}} \bigg) T_{w_\lambda} y_{\lambda'}.
\end{align*}

Now suppose that $Y \sim_r \bar{V}$ and that $Y(c,b) = Y(c+1,b)$ for some $b$. Then $1_{Y} \in \mathscr{D}_\lambda$ so that $T_{1_Y}T_{w_\lambda} = T_{1_Y w_\lambda}$.  Clearly $\mf{t}^\nu 1_Y w_\lambda$ contains two entries from the same column of $\mf{t}_\lambda$ in row $c+1$.  So $(i) 1_Y w_\lambda$ and $(j) 1_Y w_\lambda$ are in the same column of $\mf{t}_\lambda$ for some $i,j$ which lie in row $c+1$ of $\mf{t}^\nu$.  

Now, $i,j \leq \mu_1+\ldots+\mu_{c+1}$, so that any permutation which occurs in the sum $\bar{h} T_w T_{1_Y w_\lambda}$ still sends $i$ to $(i) 1_Y w_\lambda$ and $j$ to $(j) 1_Y w_\lambda$.  
Hence $x_\nu \bar{h} T_w T_{1_Y} T_{w_\lambda} = \sum_{d \in \mathscr{D}^\nu} f(d) x_\nu T_d$ for some $f(d) \in F$, where every permutation $d$ that occurs in this sum has the property that $\mf{t}^\nu d$ contains two numbers from the same column of $\mf{t}_\lambda$ in the same row. Therefore $x_\nu \bar{h} T_w T_{1_Y} T_{w_\lambda} y_{\lambda'} =0$ by Lemma \ref{nicelemma2}.  

For $1 \leq x \leq \mu_{c}$, let $V_x \sim_r \bar{V}$ be defined by 
\[V_x(a,b)= \begin{cases}  a+1 & \text{if } a=c,c+1 \text{ and } b=x,  \\
a & \text{otherwise}. \end{cases} \]
Then  
\[\hat{\Theta}_{V}(x_\lambda T_{w_\lambda} y_{\lambda'}) = x_\nu \bar{h} T_w \Bigg( \sum_{x=1}^{\mu_{c}} T_{1_{V_x}} \Bigg) T_{w_\lambda} y_{\lambda'}. \]

For  $1 \leq x \leq \mu_c$, let $U_x \sim_r \bar{U}$ be defined by 
\[U_x(a,b)= \begin{cases}  c+2 &\text{if } (a,b)=(c,x), \\
a & \text{otherwise}. \end{cases} \] 
Then  
\[\hat{\Theta}_{U}(x_\lambda T_{w_\lambda} y_{\lambda'}) = x_\nu \bar{h} T_w \Bigg( \sum_{x=1}^{\mu_{c}} T_{1_{U_x}} \Bigg) T_{w_\lambda} y_{\lambda'}. \]
We will show that, for $1 \leq x \leq \mu_{c}$, 
\[ x_\nu \bar{h} T_w T_{1_{V_x}} T_{w_\lambda} y_{\lambda'} = - x_\nu \bar{h} T_w T_{1_{U_x}} T_{w_\lambda} y_{\lambda'},\]
completing the proof of Lemma \ref{littleswap}.

Choose $x$ with $1 \leq x \leq \mu_{c}$.  Then $1_{V_x}\in \mathscr{D}_\lambda$, so $T_{1_{V_x}} T_{w_\lambda} = T_{1_{V_x} w_\lambda}$.  
Consider $\mf{t}^\nu 1_{V_x} w_\lambda$.  It is not row standard; the first entry in row $c+1$ must be moved $x-1$ places to the right.  
Let $\mf{v}_x$ be the row standard $\nu$-tableau obtained by reordering the rows of $\mf{t}^\nu 1_{V_x} w_\lambda$ and define $v_x \in \sym{n}$ by $\mf{t}^\nu v_x =\mf{v}_x$.  Then
\[ x_\nu \bar{h} T_w T_{1_{V_x}} T_{w_\lambda} y_{\lambda'} = q^{x-1}  x_\nu \bar{h} T_w T_{v_x} y_{\lambda'}.\]

Now consider  $x_\nu \bar{h} T_w T_{1_{U_x}} T_{w_\lambda} y_{\lambda'}$.  We first look at $T_{1_{U_x}} T_{w_\lambda}$.  Observe that  
\[1_{U_x} = s_{\mu_1+\ldots+\mu_{c+1}}s_{\mu_1+\ldots+\mu_{c+1}-1}\ldots s_{\mu_1+\ldots+\mu_{c-1}+x+1}.\]
For $1 \leq i \leq x-1$, let $p(i) \in \sym{n}$ be given by 
\[p(i) = s_{\mu_1+\ldots+\mu_{c+1}}s_{\mu_1+\ldots+\mu_{c+1}-1}\ldots s_{\mu_1+\ldots+\mu_{c}+i+1}  s_{\mu_1+\ldots+\mu_{c}+i-1} \ldots s_{\mu_1+\ldots+\mu_{c-1}+x+1}.\]
Then repeated application of Equation \ref{length} shows that
\[T_{1_{U_x}} T_{w_\lambda} = q^{x-1}T_{1_{U_x} w_\lambda} 
+\sum_{i=1}^{x-1} q^{i-1}(q-1) T_{p(i)w_\lambda}.\]
However, for $1 \leq i \leq x-1$, the tableau $\mf{t}^\nu p(i)w_\lambda$ contains, in row $c+1$, two entries from the same row of $\mf{t}_\lambda$.  An argument similar to that given above shows that $x_\nu \bar{h} T_w T_{p(i)w_\lambda} y_{\lambda'} =0$.
Therefore 
\[ x_\nu \bar{h} T_w T_{1_{U_x}} T_{w_\lambda} y_{\lambda'} = q^{x-1}  x_\nu \bar{h} T_w T_{1_{U_x} w_\lambda} y_{\lambda'}.\]  
It remains for the reader to convince themselves that 
$v_x s_r = 1_{U_x} w_\lambda$ for some transposition $s_r \in \sym{\lambda'}$.    Recall that $T_z y_{\lambda'} = - T_{z s_i} y_{\lambda'}$ for all $z \in \sym{n}$ and $s_i \in \sym{\lambda'}$.  Therefore  
\begin{align*}
x_\nu \bar{h} T_w T_{1_{U_x}} T_{w_\lambda} y_{\lambda'} 
& = q^{x-1}  x_\nu \bar{h} T_w T_{1_{U_x} w_\lambda} y_{\lambda'} \\
& = q^{x-1}  x_\nu \bar{h} T_w T_{v_x s_r} y_{\lambda'} \\
& = - q^{x-1}  x_\nu \bar{h} T_w T_{v_x} y_{\lambda'}\\
&= - x_\nu \bar{h} T_w T_{1_{V_x}} T_{w_\lambda} y_{\lambda'}. 
\end{align*}
\end{proof}

\begin{lemma}\label{swaprow2}
Let $A=(\mu:i_1,i_2,\ldots,i_s) \in \mathcal{T}_0(\lambda,\mu)$ be such that $i_{d+1}=d+1$ and $i_d \neq d$. 
Let $U = (\nu:i_1,i_2,\ldots,i_{d-1},d,i_{d},i_{d+2},\ldots,i_s) \in \mathcal{T}(\lambda,\nu)$.

Suppose that $\mu_{d+1}=\mu_{d+l}>\mu_{d+l+1}$ and that $i_d <i_{d+l}$.  Then $\{i_d,i_{d+2},i_{d+3},\ldots,i_{d+l-1}\} = \{d+2,d+3,\ldots,d+l\}$.
Let $V = (\nu:i_1,i_2,\ldots,i_{d-1},d,d+2,\ldots,d+l,i_{d+l},i_{d+l+1},\ldots,i_s) \in \mathcal{T}(\lambda,\nu)$.  
Then $V$ is semistandard and $\hat{\Theta}_{U}=(-1)^{i_d-d}\hat{\Theta}_{V}$.
\end{lemma}
\begin{proof}
The proof follows from Lemma \ref{littleswap}.
\end{proof}

\begin{lemma}\label{swaprow1}
Let $A=(\mu:i_1,i_2,\ldots,i_s) \in \mathcal{T}_0(\lambda,\mu)$ be such that $i_{d+1}=d+1$ and $i_d \neq d$. 
Let $U = (\nu:i_1,i_2,\ldots,i_{d-1},d,i_{d},i_{d+2},\ldots,i_s) \in \mathcal{T}(\lambda,\nu)$.

Suppose that $\mu_{d+1}=\mu_{d+l}>\mu_{d+l+1}$ and that $i_d >i_{d+l}$.  Then $\{i_{d+2},i_{d+3},\ldots,i_{d+l}\} = \{d+2,d+3,\ldots,d+l\}$.
Let $V = (\nu:i_1,i_2,\ldots,i_{d-1},d,d+2,\ldots,d+l,i_d,i_{d+l+1},\ldots,i_s) \in \mathcal{T}(\lambda,\nu)$.  
Then $V$ is semistandard and $\hat{\Theta}_{U}=(-1)^{l-1}\hat{\Theta}_{V}$.
\end{lemma}
\begin{proof}
The proof follows from Lemma \ref{littleswap}.
\end{proof}

\subsection{Conditions}
We now consider the maps $\psi_d \hat{\Theta}_A:S^\lambda \rightarrow M^\nu$.  We write $\psi_d \hat{\Theta}_A$ in terms of homomorphisms indexed by semistandard $\lambda$-tableaux of type $\nu$.  We stress that the hard work has already been done, and it is now just a question of collecting together our results.
We examine five separate cases.

\renewcommand{\thecase}{\Alph{case}}
\begin{case}[$d=s$]
Let $A=(\mu:i_1,i_2,\ldots,i_{s})\in\mathcal{T}_0(\lambda,\mu)$. Recall that $r(a')$ is defined to be such that $i_{r(a')}=a'$.
\begin{lemma} \label{startlemma}
Suppose that $r(s+1)<s$; therefore $i_s=s$.  Let \\$S=(\nu:i_1,\ldots,i_{r(s+1)-1},s,i_{r(s+1)+1},\ldots,i_{s-1},s)\in\mathcal{T}(\lambda,\nu)$.  Then $\psi_d \hat{\Theta}_A =q^{\mu_s} \hat{\Theta}_S$, and $S$ is  semistandard.
\end{lemma}
\begin{proof}
The proof follows from Proposition \ref{comb}.
\end{proof}

\begin{lemma}
Suppose that $r(s+1) =s$.  Let $S=(\nu:i_1,i_2,\ldots,i_{s-1},s)\in\mathcal{T}(\lambda,\nu)$.  Then $\psi_d \hat{\Theta}_A = [\mu_s] \hat{\Theta}_S$.  If $\mu_{s-1}=\mu_s$ and $i_{s-1}=s$ then $\hat{\Theta}_S=0$; otherwise $S$ is semistandard.  
\end{lemma}
\begin{proof}
The proof follows from Proposition \ref{comb} and Lemma \ref{kill2}. 
\end{proof}
\end{case}

\begin{case}[$2\leq d <s$ and $\mu_{d+1}=1$]Let $A=(\mu:i_1,i_2,\ldots,i_{s})\in\mathcal{T}_0(\lambda,\mu)$.
\begin{lemma}
Suppose that $r(d+1)<d$; therefore $r(d)=d$.  Let \\$S=(\nu:i_1,\ldots,i_{r(d+1)-1}, d,i_{r(d+1)+1},\ldots,i_{s})\in\mathcal{T}(\lambda,\nu)$.  Then $\psi_d \hat{\Theta}_A =q^{\mu_d} \hat{\Theta}_S$, and $S$ is  semistandard.
\end{lemma}
\begin{proof}
The proof follows from Proposition \ref{comb}.
\end{proof}

\begin{lemma}
Suppose that $r(d+1) =d$.  Let $S=(\nu:i_1\ldots,i_{d-1},d,i_{d+1},\ldots,i_{s})\in\mathcal{T}(\lambda,\nu)$.  Then $\psi_d \hat{\Theta}_A = [\mu_d] \hat{\Theta}_S$.  If $\mu_{d-1}=\mu_d$ and $i_{d-1}=d$ then $\hat{\Theta}_S=0$; otherwise $S$ is semistandard.  
\end{lemma}
\begin{proof}
The proof follows from Proposition \ref{comb} and Lemma \ref{kill2}. 
\end{proof}

\begin{lemma}
Suppose that $r(d+1)=d+1$.  If $i_d=d$ then $\psi_d \hat{\Theta}_A=0$, and if $\mu_d=\mu_{d+1}$ then $i_d=d$.  Else let $S = (i_1,\ldots,i_{d-1},d,d+2,\ldots,s+1)\in\mathcal{T}(\lambda,\nu)$.  Then $\psi_d \hat{\Theta}_A =(-1)^{i_d-d-1}\hat{\Theta}_{S}$, and $S$ is  semistandard.
\end{lemma}

\begin{proof}
The proof follows from Proposition \ref{comb} and Lemmas \ref{kill1}, \ref{swap2}, \ref{swaprow2} and \ref{swaprow1}. 
\end{proof}
\end{case}

\begin{case}[$2 \leq d<s$ and $\mu_{d+1}>1$]Let $A=(\mu:i_1,i_2,\ldots,i_{s})\in\mathcal{T}_0(\lambda,\mu)$.  Define $l$ such that $\mu_{d+1}=\mu_{d+l}>\mu_{d+l+1}$.
\begin{lemma}
Suppose that $r(d+1)<d$; therefore $r(d)=d$.  Let \\$S=(\nu:i_1,\ldots,i_{r(d+1)-1}, d,i_{r(d+1)+1},\ldots,i_{s})\in\mathcal{T}(\lambda,\nu)$.  Then $\psi_d \hat{\Theta}_A =q^{\mu_d} \hat{\Theta}_S$, and $S$ is  semistandard.
\end{lemma}
\begin{proof}
The proof follows from Proposition \ref{comb}.
\end{proof}

\begin{lemma}
Suppose that $r(d+1) =d$.  Let $S=(\nu:i_1\ldots,i_{d-1},d,i_{d+1},\ldots,i_{s})\in\mathcal{T}(\lambda,\nu)$.  
If $\mu_{d-1}=\mu_d$ and $i_{d-1}=d$ then $\psi_d \hat{\Theta}_A = 0$.   
Else $\psi_d \hat{\Theta}_A = ([\mu_d]-[\mu_{d+1}-1]) \hat{\Theta}_S = q^{\mu_{d+1}-1}[\mu_d-\mu_{d+1}+1]\hat{\Theta}_S$, and $S$ is  semistandard.
\end{lemma}
\begin{proof}
The proof follows from Proposition \ref{comb} and Lemmas \ref{swap1} and \ref{kill2}. 
\end{proof}

\begin{lemma}
Suppose that $r(d+1)=d+1$.  If $i_d=d$ then $\psi_d \hat{\Theta}_A=0$, and if $\mu_d=\mu_{d+1}$ then $i_d=d$.  
Else if $i_d<i_{d+l}$, let $S = (i_1,\ldots,i_{d-1},d,d+2,\ldots,d+l,i_{d+l},\ldots,i_{s})\in\mathcal{T}(\lambda,\nu)$;  
then $\psi_d \hat{\Theta}_A =(-1)^{i_d-d-1}q^{\mu_{d+1}-1}\hat{\Theta}_{S}$, and $S$ is  semistandard. If $i_d>i_{d+l}$, let $S = (i_1,\ldots,i_{d-1},d,d+2,\ldots,d+l,i_{d},\ldots,i_{s})\in\mathcal{T}(\lambda,\nu)$; then $\psi_d \hat{\Theta}_A =(-1)^{l}q^{\mu_{d+1}-1}\hat{\Theta}_{S}$, and $S$ is  semistandard.
\end{lemma}

\begin{proof}
The proof follows from Proposition \ref{comb} and Lemmas \ref{kill1}, \ref{swap2}, \ref{swaprow2} and \ref{swaprow1}. 
\end{proof}
\end{case}

\begin{case}[$d=1$ and $\mu_2=1$]Let $A=(\mu:i_1,i_2,\ldots,i_{s})\in\mathcal{T}_0(\lambda,\mu)$. 

\begin{lemma}
Suppose that $r(2)=1$.  Let $S=(\nu:1,3,\ldots,s+1)\in\mathcal{T}(\lambda,\nu)$.  Then $\psi_d \hat{\Theta}_A = [\mu_1+1] \hat{\Theta}_S$, and $S$ is  semistandard.  
\end{lemma}

\begin{proof}
The proof follows from Proposition \ref{comb}. 
\end{proof}

\begin{lemma}
Suppose that $r(2)=2$. 
Let $S = (1,3,\ldots,s+1)\in\mathcal{T}(\lambda,\nu)$.  
Then $\psi_d \hat{\Theta}_A =(-1)^{i_1}\hat{\Theta}_{S}$, and $S$ is  semistandard.
\end{lemma}

\begin{proof}
The proof follows from Proposition \ref{comb} and Lemmas \ref{swap2}, \ref{swaprow2} and \ref{swaprow1}.
\end{proof}

\end{case}

\begin{case}[$d=1$ and $\mu_2>1$]Let $A=(\mu:i_1,i_2,\ldots,i_{s})\in\mathcal{T}_0(\lambda,\mu)$.  Define $l$ such that $\mu_{d+1}=\mu_{d+l}>\mu_{d+l+1}$.
\begin{lemma}
Suppose that $r(2)=1$.  Let $S=(\nu:1,i_2,\ldots,i_s)\in\mathcal{T}(\lambda,\nu)$.  
Then $\psi_d \hat{\Theta}_A = ([\mu_1+1]-[\mu_{2}-1]) \hat{\Theta}_S = q^{\mu_{2}-1}[\mu_1-\mu_{2}+2]\hat{\Theta}_S$, and $S$ is  semistandard. 
\end{lemma}
\begin{proof}
The proof follows from Proposition \ref{comb} and Lemma \ref{swap1}. 
\end{proof}

\begin{lemma}\label{endlemma}
Suppose that $r(2)=2$.
If $i_1<i_{l+1}$, let $S = (1,2,\ldots,l+1,i_{l+1},i_{l+2},\ldots,i_{s})\in\mathcal{T}(\lambda,\nu)$;  
then $\psi_d \hat{\Theta}_A =(-1)^{i_1}q^{\mu_{2}-1}\hat{\Theta}_{S}$, and $S$ is  semistandard. If $i_1>i_{l+1}$, let $S = (1,2,\ldots,l+1,i_{1},i_{l+2},\ldots,i_s)\in\mathcal{T}(\lambda,\nu)$; then $\psi_d \hat{\Theta}_A =(-1)^{l}q^{\mu_{2}-1}\hat{\Theta}_{S}$ and $S$ is  semistandard.
\end{lemma}

\begin{proof}
The proof follows from Proposition \ref{comb} and Lemmas \ref{swap2}, \ref{swaprow2} and \ref{swaprow1}.
\end{proof}
\end{case}

Proposition \ref{sumconds} summarises the results of Lemmas \ref{startlemma} to \ref{endlemma}.
 
\begin{proposition}\label{sumconds}
For $1\leq d \leq s$, $\psi_d \hat{\Theta}=0$  if and only if the following equations all hold.
\begin{itemize}
\item If $d>1$ and $\mu_d=\mu_{d+1}$:
\begin{multline*}
q^{\mu_d}\, f(\mu:j_1,\ldots,j_{\check{r}(d)-1},d+1,j_{\check{r}(d)+1},\ldots,j_{d-1},d,d+2,,\ldots,d+l,j_{d+l},\ldots,j_s)\\
+q^{\mu_{d+1}-1}\,f(\mu:j_1,\ldots,j_{\check{r}(d)-1},d,j_{\check{r}(d)+1},\ldots,j_{d-1},d+1,d+2,\ldots,d+l,j_{d+l},\ldots,j_s)=0.
\end{multline*}
\item If $d>1$ and $\mu_d>\mu_{d+1}=\mu_{d+l}>\mu_{d+l+1}$:
\begin{align*}
q^{\mu_d}\,&f(\mu:j_1,\ldots,j_{\check{r}(d)-1},d+1,j_{\check{r}(d)+1},\ldots,j_{d-1},d,d+2,\ldots,d+l,j_{d+l},\ldots,j_s) \\
+ &q^{\mu_{d+1}-1}[\mu_d-\mu_{d+1}+1]\\
& \quad \quad f(\mu:j_1,\ldots,j_{\check{r}(d)-1},d,j_{\check{r}(d)+1},\ldots,j_{d-1},d+1,d+2,\ldots,d+l,j_{d+l},\ldots,j_s) \\
- &q^{\mu_{d+1}-1}f(\mu:j_1,\ldots,j_{\check{r}(d)-1},d,j_{\check{r}(d)+1},\ldots,j_{d-1},d+2,d+1,\ldots,d+l,j_{d+l},\ldots,j_s) \\
&\ldots \\
+ &(-1)^{l-1} q^{\mu_{d+1}-1}
f(\mu:j_1,\ldots,j_{\check{r}(d)-1},d,j_{\check{r}(d)+1},\ldots,j_{d-1},d+l,d+1,\ldots,d+l-1,j_{d+l},\ldots,j_s)\\
+&(-1)^l  q^{\mu_{d+1}-1} f(\mu:j_1,\ldots,j_{\check{r}(d)-1},d,j_{\check{r}(d)+1},\ldots,j_{d-1},j_{d+l},d+1,\ldots,d+l-1,d+l,\ldots,j_s)\\&=0.
\end{align*}
\item If $d=1$ and $\mu_2=\mu_{l+1}>\mu_{l+2}$:
\begin{align*}
[\mu_1-\mu_2+2] \, &f(\mu: 2,3,4,\ldots,l+1,j_{l+1},j_{l+2},\ldots,j_s) - f(\mu: 3,2,4,\ldots,l+1,j_{l+1},j_{l+2},\ldots,j_s) \\
&+f(\mu: 4,2,3,\ldots,l+1,j_{l+1},j_{l+2},\ldots,j_s)\\
& \quad \ldots \\
&+(-1)^{l+1} f(\mu:l+1,2,3,\ldots,l,j_{l+1},j_{l+2},\ldots,j_s)\\
&+(-1)^{l} f(\mu:j_{l+1},2,3,\ldots,l,l+1,j_{l+2},\ldots,j_s)=0.
\end{align*}
\end{itemize}
\end{proposition}

\subsection{Summary}
We now write down a map which will satisfy these conditions.

\begin{definition}\label{defnmap}
Let $A \in \mathcal{T}_0(\lambda,\mu)$. For $2 \leq i \leq s$, define $A(i) \in F$ by
\[A(i) = \begin{cases} 1 &\text{if } A(i,\lambda_i) \neq i, \\
-q^{-1} &\text{if } A(i,\lambda_i) = i \text{ and } \lambda_i = \lambda_{i+1}, \\
-q^{-(\lambda_i+s-i)}\,[\lambda_i+s-i] &\text{if } A(i,\lambda_i) = i \text{ and } \lambda_i \neq \lambda_{i+1}. \end{cases}\]
Set 
\[f(A) = \prod_{i=2}^s A(i)\]
and define $\hat{\Theta}:S^\lambda \rightarrow M^\mu$ by
\[\hat{\Theta} = \sum_{A \in \mathcal{T}_0(\lambda,\mu)}f(A)\hat{\Theta}_A. \]
Note that $f(\mu:2,3,\ldots,s+1)=1$, so that $\hat{\Theta} \neq 0$.
\end{definition}

\begin{theorem}\label{inimage}
Suppose $e \mid \mu_1+s$.  Then $\Image(\hat{\Theta}) \subseteq S^\lambda$.
\end{theorem}

\begin{proof}We consider the conditions of Proposition \ref{sumconds}.
\begin{itemize}
\item Suppose $d>1$ and $\mu_d=\mu_{d+1}$.  Consider
\begin{multline}\label{sum1}
q^{\mu_d}\, f(\mu:j_1,\ldots,j_{\check{r}(d)-1},d+1,j_{\check{r}(d)+1},\ldots,j_{d-1},d,d+2,,\ldots,d+l,j_{d+l},\ldots,j_s) \\
+q^{\mu_{d+1}-1}\,f(\mu:j_1,\ldots,j_{\check{r}(d)-1},d,j_{\check{r}(d)+1},\ldots,j_{d-1},d+1,d+2,\ldots,d+l,j_{d+l},\ldots,j_s). 
\end{multline}
We may ignore all values of $A(j)$ except $j=\check{r}(d)$ and $j=d$. Hence, for some $C \in F$,
\[(\ref{sum1}) = C (q^{\mu_d} (-q^{-1})+q^{\mu_{d+1}-1})=0.\]
\item Suppose $d>1$ and $\mu_d>\mu_{d+1}$.  Then there exists $C \in F$ such that 
\begin{align*}
q&^{\mu_d}\,f(\mu:j_1,\ldots,j_{\check{r}(d)-1},d+1,j_{\check{r}(d)+1},\ldots,j_{d-1},d,d+2,\ldots,d+l,j_{d+l},\ldots,j_s) \\
&+ q^{\mu_{d+1}-1}[\mu_d-\mu_{d+1}+1]\\
&\quad \quad f(\mu:j_1,\ldots,j_{\check{r}(d)-1},d,j_{\check{r}(d)+1},\ldots,j_{d-1},d+1,d+2,\ldots,d+l,j_{d+l},\ldots,j_s) \\
&- q^{\mu_{d+1}-1}f(\mu:j_1,\ldots,j_{\check{r}(d)-1},d,j_{\check{r}(d)+1},\ldots,j_{d-1},d+2,d+1,\ldots,d+l,j_{d+l},\ldots,j_s) \\
&\quad \ldots \\
&+ (-1)^{l-1} q^{\mu_{d+1}-1}f(\mu:j_1,\ldots,j_{\check{r}(d)-1},d,j_{\check{r}(d)+1},\ldots,j_{d-1},d+l,d+1,\ldots,d+l-1,j_{d+l},\ldots,j_s)\\
&+(-1)^l  q^{\mu_{d+1}-1}f(\mu:j_1,\ldots,j_{\check{r}(d)-1},d,j_{\check{r}(d)+1},\ldots,j_{d-1},j_{d+l},d+1,\ldots,d+l-1,d+l,\ldots,j_s)\\
 = &C\Big(q^{\mu_d}(-q^{-(\mu_d+s-d)}[\mu_d+s-d]) + q^{\mu_{d+1}-1}[\mu_d-\mu_{d+1}+1] -  q^{\mu_{d+1}-1}(-q^{-1}) + \ldots + \\
& \quad (-1)^{-(l-1)} q^{\mu_{d+1}-1}((-q)^{l-1}) + (-1)^{l} q^{\mu_{d+1}-1} (-(-q)^{-(l-1)}q^{-(\mu_{d+l}+s-d-l)}[\mu_{d+l}+s-d-l])\Big)\\
=&C\Big(-q^{-(s-d)}[\mu_d+s-d] + q^{\mu_{d+1}-1}\Big([\mu_d-\mu_{d+1}+1] + q^{-(l-1)}[l-1] \\
&\quad +q^{-(\mu_{d+1}+s-d-l)}[\mu_{d+1}+s-d-l]\Big) \Big)\\ 
=&0.
\end{align*}
\item Suppose $d=1$.  Then there exists $C \in F$ such that 
\begin{align*}
[\mu_1-\mu_2+2] \, &f(\mu: 2,3,4,\ldots,l+1,j_{l+1},j_{l+2},\ldots,j_s) - f(\mu: 3,2,4,\ldots,l+1,j_{l+1},j_{l+2},\ldots,j_s) \\
&+f(\mu: 4,2,3,\ldots,l+1,j_{l+1},j_{l+2},\ldots,j_s)\\
& \quad \ldots \\
&+(-1)^{l+1} f(\mu:l+1,2,3,\ldots,l,j_{l+1},j_{l+2},\ldots,j_s)\\
&+(-1)^{l} f(\mu:j_{l+1},2,3,\ldots,l,l+1,j_{l+2},\ldots,j_s) \\
=& C\Big( [\mu_1-\mu_2+2] +q^{-(l-1)}[l-1]+q^{-(l-1)}q^{-(\mu_2+s-l-1)}[\mu_2+s-l-1]\Big)\\
=&C \Big(q^{-(\mu_2+s-2)} [\mu_1+s]\Big) \\
=&0.
\end{align*}
\end{itemize}
\end{proof}

\begin{theorem}\label{image} 
Suppose that 
\begin{align*}
\mu&=(\mu_1,\mu_2,\ldots,\mu_s,1), \\
\lambda&=(\mu_1+1,\mu_2,\ldots,\mu_s).
\end{align*}
are partitions of $n$ and that $e \neq 2$ or $\lambda$ is 2-regular.  Then
\[\dim(\Hom_\h(S^\lambda,S^\mu)) = \begin{cases} 1 & \text{if } e \mid \mu_1 +s, \\
0 & \text{otherwise}. \end{cases}\]
\end{theorem}

\begin{proof}
By Theorem \ref{semistand}, every map $\hat{\Theta}:S^\lambda \rightarrow M^\mu$ is a sum of semistandard homomorphisms.
Define a total order $\preceq$ on $\mathcal{T}_0(\lambda,\mu)$ by saying that $A_i=(\mu:i_1,i_2,\ldots,i_s) \preceq A_j=(\mu:j_1,j_2,\ldots,j_s)$ if $A_i=A_j$, or there exists $b$ such that $i_a=j_a$ for $a<b$ and $i_b<j_b$.  If $A_i \preceq A_j$ and $A_i \neq A_j$, write $A_i \prec A_j$.

Let $A=(\mu:i_1,i_2,\ldots,i_s)\in  \mathcal{T}_0(\lambda,\mu)$.  Suppose $A \neq (\mu:2,3,\ldots,s+1)$.  Then there exists $b$ with $2 \leq b \leq s$ such that $i_b =b$.  Choose $b$ maximal such that $i_b =b$; then $i_{r(b+1)} = b+1$ for some $r(b+1)<b$. Consider the coefficient of 
$(\nu:i_1,\ldots,i_{r(b+1)-1},b,i_{r(b+1)+1},\ldots,i_{b-1},b,i_{b+1},\ldots,i_s)$
in $\psi_b \hat{\Theta}$.  We get a relation 
\[ q^{\mu_b}  f(A) + \sum_{B \prec A} g(B) f(B)  =0\]
for some $g(B) \in F$.  Hence the space of homomorphisms which satisfy the conditions of Proposition \ref{sumconds} when $d>1$ is at most one dimensional.  The proof of Theorem \ref{image} shows that the map $\hat{\Theta}$ of Definition \ref{defnmap} always satisfies the conditions of Proposition \ref{sumconds} when $d>1$; and satisfies the conditions of Proposition \ref{sumconds} when $d=1$ if and only if $e \mid \mu_1+s$.

Of course, the fact that $\dim(\Hom_\h(S^\lambda,S^\mu)) = 0$ if $e \nmid \mu_1+s$ can also be deduced from the Nakayama conjecture (see~\cite{M:ULect}, Corollary 5.38).   
\end{proof}

\begin{theorem}\label{addthm}
Suppose that
\begin{align*}
\xi&=(\xi_1,\ldots,\xi_{a-1},\xi_a,\xi_{a+1},\ldots,\xi_{b-1},\xi_b,\xi_{b+1},\ldots,\xi_r), \\
\eta&=(\xi_1,\ldots,\xi_{a-1},\xi_a+1,\xi_{a+1},\ldots,\xi_{b-1},\xi_b-1,\xi_{b+1},\ldots,\xi_r)
\end{align*} are partitions of $n$.
Let $A \in \mathcal{T}_0(\eta,\xi)$. For $a< i <b$, define $A(i) \in F$ by
\[A(i) = \begin{cases} 1 &\text{if } A(i,\eta_i) \neq i, \\
-q^{-1} &\text{if } A(i,\eta_i) = i \text{ and } \eta_i = \eta_{i+1}, \\
-q^{-(\eta_i-\eta_b+b-i-1)}\,[\eta_i-\eta_b+b-i-1] &\text{if } A(i,\eta_i) = i \text{ and } \eta_i \neq \eta_{i+1}. \end{cases}\]
Set 
\[f(A) = \prod_{i=a+1}^{b-1} A(i)\]
and define $0 \neq \hat{\Theta}:S^\eta \rightarrow M^\xi$ by
\[\hat{\Theta} = \sum_{A \in \mathcal{T}_0(\eta,\xi)}f(A)\hat{\Theta}_A. \]
Suppose $e \mid \xi_a-\xi_b+b-a+1$.
Then $\Image(\hat{\Theta}) \subseteq S^\xi$.
\end{theorem}
\begin{proof}
The proof follows along the lines of the proof of Theorem \ref{inimage}.
\end{proof}

\begin{corollary}\label{addback}Take $\xi$ and $\eta$ as in Theorem \ref{addthm} and suppose that $e \neq 2$ or $\eta$ is 2-regular.  
Then
\[\dim(\Hom_\h(S^\eta,S^\xi)) = \begin{cases} 1 & \text{if } e \mid \xi_a-\xi_b+b-a+1, \\
0 & \text{otherwise}. \end{cases}\]
\end{corollary}

\begin{corollary}Take $\xi$ and $\eta$ as in Theorem \ref{addthm} and suppose that $e = 2$ and $\eta$ is not 2-regular.  
Then
\[\dim(\Hom_\h(S^\eta,S^\xi)) = \begin{cases} c \geq 1 & \text{if } e \mid \xi_a-\xi_b+b-a+1, \\
0 & \text{otherwise}. \end{cases}\]
\end{corollary}

\numberwithin{equation}{section}
\section{Reducible Specht modules}
The reducible Specht modules for the symmetric group algebras have been classified in the series of papers~\cite{James:CartC,Fayers:red,Fayers:irred,JM:Specht2,Lyle:irred}.  We are now in a position to complete the classification of the reducible Specht modules for the Hecke algebra $\h_{F,q}(\sym{n})$ when $e \neq 2$, verifying the conjecture of James and Mathas~\cite{M:ULect}, Conjecture~5.47.  

Let $\lambda$ be a partition and recall that the diagram of $\lambda$ is the set of nodes
\[[\lambda] = \{(i,j)\mid 1 \leq i \mbox{ and } 1 \leq j\leq \lambda_i\}.\]
For each node $(i,j)$ in $[\lambda]$, we define the $(i,j)$-hook length 
$h^\lambda_{ij}=\lambda_i-i+\lambda_j'-j+1$.
Define $\nu_{e,p}:\N\rightarrow\Z$ by
$$\nu_{e,p}(h)=\begin{cases}
        \nu_p(\tfrac he)+1 &\text{if $e$ divides $h$},\\
              0 &\text{otherwise,}
\end{cases}$$
where $\nu_p(k)$ is maximal such that $p^{\nu_p(k)} \mid k$. If $p=0$ then set $\nu_p(k)=0$, for all $k$.

\begin{definition} \label{epdefn} 
A partition $\lambda$ is said to be $(e,p)$-reducible if there exist
nodes $(a,i)$, $(a,j)$ and $(b,i)$ in $[\lambda]$ such that
$\nu_{e,p}(h^\lambda_{ai})>0$, and
$\nu_{e,p}(h^\lambda_{aj})\ne\nu_{e,p}(h^\lambda_{ai})
                     \ne\nu_{e,p}(h^\lambda_{bi})$.
\end{definition}

We begin by describing some particular reducible Specht modules.

\begin{theorem}\label{irredone}
Suppose that $e \neq 2$.  The Specht module $S^{\lambda}$ is reducible if there exist nodes $(a,i)$, $(a,j)$ and $(b,i)$ in $[\lambda]$ such that
$e \mid (h^\lambda_{ai})$, and
$e \nmid (h^\lambda_{aj})$ and 
$ e \nmid (h^\lambda_{bi})$.
\end{theorem}

\begin{proof}
Theorem \ref{irredone} was initially proved for the symmetric group algebras~\cite{Lyle:irred}, Theorem 2.16.  It used a result of Brundan and Kleshchev~\cite{kl}, Theorem 2.13 which was originally stated for the symmetric group algebras, however the proof given in~\cite{kl} also works for arbitrary Hecke algebras.  We are grateful to Alexander Kleshchev for this information.  Furthermore, the proof of~\cite{Lyle:irred}, Theorem 2.16 relied on the existence of non--zero homomorphisms between certain Specht modules.  These Specht modules were indexed by partitions which fulfilled the conditions of Theorem \ref{mainextra}.  Given the main result of this paper, Theorem \ref{mainextra}, the proof of Theorem \ref{irredone} follows immediately from the corresponding proof in~\cite{Lyle:irred}.  
\end{proof}

\begin{theorem}\label{irred}
Suppose that $e \neq 2$.  The Specht module $S^{\lambda}$ is reducible if and only if $\lambda$ is $(e,p)$-reducible.
\end{theorem}  

\begin{proof}
A proof that if $\lambda$ is not  $(e,p)$-reducible then $S^\lambda$ is irreducible is given by Fayers~\cite{Fayers:irred} in the cases that $q=1$ or $F$ is a field of characteristic zero. This proof has been generalised to arbitrary Hecke algebras in~\cite{JLM:Rouquier}.  
Combining Theorem \ref{irredone} with the results of~\cite{Fayers:red} shows that if the partition $\lambda$ is $(e,p)$-reducible then the Specht module $S^\lambda$ is reducible, completing the proof of Theorem \ref{irred}.
\end{proof}

\end{document}